%% file: Nonlinear-thermomechanical.tex
\documentclass[sort&compress,11pt]{elsarticle}
\usepackage{verbatim,amsmath,amssymb}
\usepackage{epsfig,float,color}
\usepackage{geometry}
\usepackage{setspace}
\usepackage{natbib}
\usepackage{amsthm,mathrsfs,amsfonts,dsfont}
\usepackage{hyperref}
\usepackage[yyyymmdd,hhmmss]{datetime}
\usepackage[capitalise]{cleveref}
\usepackage{caption}
\usepackage{subcaption}
\usepackage{tablefootnote}
\usepackage[section]{placeins}
\usepackage{graphicx}
\usepackage{fontenc}
\usepackage{mmap}
\usepackage{accents}
\usepackage{url}

\definecolor{darkblue}{rgb}{0,0,1}
\hypersetup{pdftex=true, colorlinks=true, breaklinks=true, linkcolor=darkblue, menucolor=darkblue, pagecolor=darkblue, citecolor=darkblue, urlcolor=darkblue}

\input{neco_updated.tex}

\graphicspath{ {/} }

\begin{document}

\begin{center}
\Large{\bf{A nonlinear thermomechanical formulation for anisotropic volume and surface continua
}}\\

\end{center}

\begin{center}

\large{Reza Ghaffari\footnote{Corresponding author, email: ghaffari@aices.rwth-aachen.de} and Roger A. Sauer\footnote{Email: sauer@aices.rwth-aachen.de}}\\
\vspace{4mm}

\small{\textit{
Aachen Institute for Advanced Study in Computational Engineering Science (AICES), \\
RWTH Aachen University, Templergraben 55, 52056 Aachen, Germany}}
%
%



\end{center}

\vspace{3mm}


\rule{\linewidth}{.15mm}
{\bf Abstract:}
A thermomechanical, polar continuum formulation under finite strains is proposed for anisotropic materials using a multiplicative decomposition of the deformation gradient. First, the kinematics and conservation laws for three dimensional, polar and non-polar continua are obtained. Next, these kinematics are connected to their corresponding counterparts for surface continua based on Kirchhoff-Love kinematics. Likewise, the conservation laws for Kirchhoff-Love shells are derived from their three dimensional counterparts. From this, the weak forms are obtained for three dimensional non-polar continua and Kirchhoff-Love shells. These formulations are expressed in tensorial form so that they can be used in both curvilinear and Cartesian coordinates. They can be used to model anisotropic crystals and soft biological materials, and they can be extended to other field equations, like Maxwell's equations to model thermo-electro-magneto-mechanical materials.
\rule{\linewidth}{.15mm}
\section*{List of important symbols}
\vspace{-2mm}
\begin{tabbing}
$\bone$ \qquad~~~~~~~~~ \=  identity tensor in $\bbR^3$  \\
$\bA_\alpha$  \> co-variant tangent vectors of $\mcalS_{0}$ at point $\bXz$ ; $\alpha=1,2$ \\
$\ba_\alpha$  \> co-variant tangent vectors of $\mcalS$ at point $\bxz$ ; $\alpha=1,2$ \\
$\baTz_{\alpha}$  \> co-variant tangent vectors of $\mcalSTz$ at point $\bXTz$ ; $\alpha=1,2$ \\
$\bA^\alpha$ \> contra-variant tangent vectors of $\mcalS_{0}$ at point $\bXz$ ; $\alpha=1,2$\\
$\ba^\alpha$  \> contra-variant tangent vectors of $\mcalS$ at point $\bxz$ ; $\alpha=1,2$  \\
$\baTz^{\alpha}$  \> contra-variant tangent vectors of $\mcalSTz$ at point $\bXTz$ ; $\alpha=1,2$ \\
$A_{\alpha\beta}$ \> co-variant metric of $\mcalS_{0}$ at point $\bXz$ ; $\alpha,\beta=1,2$\\
$a_{\alpha\beta}$ \> co-variant metric of $\mcalS$ at point $\bxz$ ; $\alpha,\beta=1,2$\\
$\aTz_{\alpha\beta}$ \> co-variant metric of $\mcalS_{\text{T}}$  at point $\bXTz$ ; $\alpha,\beta=1,2$\\
$A^{\alpha\beta}$ \> contra-variant metric of $\mcalS_{0}$ at point $\bXz$ ; $\alpha,\beta=1,2$\\
$a^{\alpha\beta}$ \> contra-variant metric of $\mcalS$ at point $\bxz$ ; $\alpha,\beta=1,2$\\
$\mcalB_{0}$, $\mcalBTd$, $\mcalB$ \>  reference, intermediate, current configuration of three dimensional continua\\
$b_{0\,\alpha\beta}$ \>  co-variant curvature tensor components of $\mcalS_{0}$ at point $\bXz$ ; $\alpha,\beta=1,2$\\
$\buab$ \> co-variant curvature tensor components of $\mcalS$ at point $\bxz$ ; $\alpha,\beta=1,2$\\
$b_{0}^{\alpha\beta}$ \>  contra-variant curvature tensor components of $\mcalS_{0}$ at point $\bXz$ ; $\alpha,\beta=1,2$\\
$\bab$ \> contra-variant curvature tensor components of $\mcalS$ at point $\bxz$ ; $\alpha,\beta=1,2$\\
$\bb_{0}$ \> curvature tensor of $\mcalS_{0}$ at point $\bXz$ \\
$\bb$ \> curvature tensor of $\mcalS$ at point $\bxz$ \\
$\bbTz$ \> curvature tensor of $\mcalSTz$ at point $\bXTz$ \\
$\bb^{\flat\lhd}$ \> pull back of $\bb$, i.e.~$\bb^{\flat\lhd}=\bFz^{\text{T}}\,\bb\,\bFz$\\
$\bbTz^{\flat\lhd(\bFTz)}$ \> pull back of $\bbTz$, i.e.~$\bbTz^{\flat\lhd(\bFTz)}=\bFTz^{\text{T}}\,\bbTz\,\bFTz$\\
$\bc$ \> body force couple (per unit mass)\\
$\bCz$ \>  surface right Cauchy–Green tensor \\
$\bCd$, $\bCd_{\text{e}}$, $\bCd_{\text{T}}$ \>  right Cauchy–Green tensor and its elastic and thermal part \\
$\bCcaz$, $\bCcaez$, $\bCcaTz$ \> surface right Cauchy–Green tensor and its elastic and thermal part \\
$\bCaz$, $\bCaez$, $\bCaTz$ \> in-plane surface right Cauchy–Green tensor and its elastic and thermal part \\
$\bbCd$ \> elasticity tensor\\
$\bbCd^{\text{L}}$ \> elasticity tensor (with rearranged order of components)\\
$\bdd$ \> rate of deformation\\
$\GammaId{i}{j}{k}$\,($\GammaIz{\alpha}{\beta}{\gamma}$) \> Christoffel symbols of the second kind of $\mcalB_{0}$ ($\mcalS_{0}$) ; $i,j,k=1,2,3$ ($\alpha,\beta=1,2$)\\
$\Gammad{i}{j}{k}$\,($\Gammaz{\alpha}{\beta}{\gamma}$) \> Christoffel symbols of the second kind of $\mcalB$ ($\mcalS$) ; $i,j,k=1,2,3$ ($\alpha,\beta=1,2$)\\
$\dif A$ ($\dif a$) \> differential surface element on $\mcalS_{0}$ ($\mcalS$)\\
$\dif S$ ($\dif s$)\> differential surface element on $\mcalB_{0}$ ($\mcalB$)\\
$\dif m$ \> differential mass element \\
$\dif V$ ($\dif v$) \> differential volume element in $\mcalB_{0}$ ($\mcalB$)\\
$\dif l$ \> differential line element on $\mcalS$\\
$\delta \bullet$ \> variation of $\bullet$\\
$\delta_{ij}, \delta^i_{j}, \delta_i^j$ \> Kronecker delta ; $i,j=1,2,3$\\
$\bed^{(n)}$ \> generalized Eulerian strain\\
$\bEd^{(n)}$ \> generalized Lagrangian strain\\
$\bsEd$ \> permutation tensor\\
$\bff$ \> body force (per unit mass)\\
$\bFd$, $\bFed$, $\bFTd$ \> deformation gradient and its elastic and thermal part\\
$\bFz$ \> in-plane surface deformation gradient of $\mcalS$, i.e.~$\bFz=\bFaz(\xi=0)$\\
$\bFIz$ \> surface deformation gradient of $\mcalS$, i.e.~$\bFIz=\bFcaz(\xi=0)$\\
$\bFaz$, $\bFaez$, $\bFaTz$ \> in-plane surface deformation gradient of $\mcalScaz$ and its elastic and thermal part\\
$\bFcaz$, $\bFcaez$, $\bFcaTz$ \> surface deformation gradient of $\mcalScaz$ and its elastic and thermal part\\
$\bG_i$ \> co-variant tangent vectors of $\mcalB_{0}$ at point $\bXd$; $i=1,2,3$ \\
$\bg_i$ \> co-variant tangent vectors of $\mcalB$ at point $\bxd$; $i=1,2,3$ \\
$\bgTdui{i}$ \> co-variant tangent vectors of $\mcalB_\text{T}$ at point $\bXTd$; $i=1,2,3$ \\
$\bGcaz_\alpha$  \> co-variant tangent vectors of $\mcalScaz_{0}$  at point $\bXcaz$; $\alpha=1,2$\\
$\bgcaz_\alpha$ \> co-variant tangent vectors of $\mcalScaz$ at point $\bxcaz$; $\alpha=1,2$\\
$\bgcaTz_{\alpha}$ \> co-variant tangent vectors of $\mcalScaz_{\text{T}}$ at point $\bXcaTz$; $\alpha=1,2$\\
$\bG^i$\> contra-variant tangent vectors of $\mcalB_{0}$ at point $\bXd$; $i=1,2,3$\\
$\bg^i$ \> contra-variant tangent vectors of $\mcalB$ at point $\bxd$; $i=1,2,3$ \\
$\bgTdi{i}$ \> contra-variant tangent vectors of $\mcalB_\text{T}$ at point $\bXTd$; $i=1,2,3$ \\
$\bGcaz^\alpha$ \> contra-variant tangent vectors of $\mcalScaz_{0}$ at point $\bXcaz$; $\alpha=1,2$ \\
$\bgcaz^\alpha$ \> contra-variant tangent vectors of $\mcalScaz$ at point $\bxcaz$; $\alpha=1,2$ \\
$\bgcaTz^{\alpha}$ \> contra-variant tangent vectors of $\mcalScaz_{\text{T}}$ at point $\bXcaTz$; $\alpha=1,2$\\
$G_{ij}$ \> co-variant metric of $\mcalB_{0}$ at point $\bXd$; $i,j=1,2,3$\\
$g_{ij}$\> co-variant metric of $\mcalB$ at point $\bxd$; $i,j=1,2,3$ \\
$\gTdui{ij}$\> co-variant metric of $\mcalBTd$ at point $\bXTd$; $i,j=1,2,3$\\
$\Gcaz_{\alpha\beta}$ \> co-variant metric of $\mcalScaz_{0}$ at point $\bXcaz$; $\alpha,\beta=1,2$\\
$\gcaz_{\alpha\beta}$ \> co-variant metric of $\mcalScaz$ at point $\bxcaz$; $\alpha,\beta=1,2$\\
$G^{ij}$ \> contra-variant metric of $\mcalB_{0}$ at point $\bXd$; $i,j=1,2,3$\\
$g^{ij}$ \> contra-variant metric of $\mcalB$ at point $\bxd$; $i,j=1,2,3$ \\
$H_{0}$ \> mean curvature of $\mcalS_{0}$ at $\bXz$\\
$H$ \> mean curvature of $\mcalS$ at $\bxz$\\
$\bHd$ \> $:=\bFTTd^{-\text{T}}\,\bCd\,\bFTd^{-\text{T}}$\\
$\bi$ \> surface identity tensor on $\mcalS$\\
$\bI$ \> surface identity tensor on $\mcalS_{0}$\\
$\Jz$ \> area change between $\mcalS_{0}$ and $\mcalS$\\
$\JTz$ \> area change between $\mcalS_{0}$ and $\mcalSTz$\\
$\Jd$ \> volume change between $\mcalB_{0}$ and $\mcalB$\\
$\JTd$ \> volume change between $\mcalB_{0}$ and $\mcalBTd$, i.e.~$\JTd=\detd\,\bFTd$\\
$\bkd$ \> conductivity tensor\\
$\kappa_{0}$ ($\kappa$) \> Gaussian curvature of $\mcalS_{0}$ ($\mcalS$) at $\bXz$ ($\bxz$)\\
$\kappa_{0\,1}$, $\kappa_{0\,2}$ \> principal curvatures of $\mcalS_{0}$ at $\bXz$\\
$\kappa_{1}$, $\kappa_{2}$ \> principal curvatures of $\mcalS$ at $\bxz$\\
$\bkappa$ \> curvature change relative to intermediate configuration\\
$\kappa_{\alpha\beta}$ \> co-variant components of curvature change relative to intermediate configuration\\
$\sK$ \> total kinetic energy\\
$\blcaz$ \> rate of deformation of $\mcalScaz$\\
$\blIz$  \> rate of deformation of $\mcalS$, i.e.~$\blIz=\blcaz(\xi=0)$ \\
$\bld$, $\bled$, $\blTd$ \> rate of deformation of $\mcalB$ and its elastic and thermal part\\
$\lambda_1$, $\lambda_2$\> principal surface stretches of $\mcalS$ at $\bx$\\
$\lambda_3$, $\lambda_{3\,\text{e}}$, $\lambda_{3\,\text{T}}$ \> stretch along $\bn$ and its elastic and thermal part\\
$m_{\nu}$, $m_{\tau}$ \> bending moment components acting at $\bxz \in \partial\mcalS$\\
$\Mab$ \>contra-variant bending moment components\\
$\Mab_{0}$ \> $=J\,\Mab$\\
$\bmz$ \> surface moment vector acting at $\bxz\in\partial\mcalS$\\
$\bmd$ \>  moment vector acting at $\bxd\in\partial\mcalB$\\
$\bmuz$ \> surface moment tensor\\
$\bmurz$ \> rotated surface moment tensor, i.e.~$\bmurz^{\text{T}}=\bn\times\bmuz^{\text{T}}$\\
$\bmurd$ \> moment tensor\\
$\bmurId$ \> two-point form of the moment tensor, i.e.~$\bmurId=\Jd\,\bFz^{-1}\,\bmurd$\\
$\bmuz^{\lhd\sharp}$ \> pull back of the surface moment tensor, i.e.~$\bmuz^{\lhd\sharp}=\bFz^{-1}\,\bmuz\bFz^{-\text{T}}$\\
$\Nab$ \> total, contra-variant, in-plane stress components\\
$\bN$ \> surface normal of $\mcalS_{0}$ at $\bXz$\\
$\bn$ \> surface normal of $\mcalS$ at $\bxz$\\
$\bn_{\text{T}}$ \> surface normal of $\mcalSTz$ at $\bXTz$\\
$\bolds{\sV}$ \> normal vector on $\partial\mcalS_{0}$ or $\partial\mcalB_{0}$\\
$\bnu$ \> normal vector on $\partial\mcalS$  or $\partial\mcalB$\\
$\bvd$ \> velocity vector \\
$\btau_{\!\text{v}}$ \> tangent vector on $\partial\mcalS$, i.e.~$\btau_{\!\text{v}}=\bn\times\bnu$ \\
$\ud$ ($\uz$) \> internal energy per unit mass, functional for $\mcalB$ ($\mcalS$)\\
$\sU$ \> total internal energy\\
$\xi^{i}$ ($\xi^{\alpha}$) \> convective coordinates on $\mcalB$ ($\mcalS$); $i=1,2,3$ ($\alpha=1, 2$)\\
$\xi$ \> thickness coordinate\\
$\psid$ ($\psiz$) \>Helmholtz free energy per unit mass, functional for $\mcalB$ ($\mcalS$)\\
$\bPd$ \> first Piola-Kirchhoff stress tensor\\
$\sP_{\text{ext}}$ \> total rate of external mechanical power\\
$\bqd$\> heat flux vector \\
$\bQd$\> pull back of heat flux vector, i.e.~$\bQd=\bqd\,\Jd\,\bFd^{-\text{T}}$ \\
$\sQ$ \> total rate of heat input\\
$r$ \> heat source per unit mass\\
$\rhoId$ ($\rhoIz$) \> density $\mcalB_{0}$ ($\mcalS_{0}$) at $\bXd$ ($\bXz$)\\
$\rhod$ ($\rhoz$) \> density $\mcalB$ ($\mcalS$) at $\bxd$ ($\bxz$)\\
$\sz$ \> entropy per unit of mass, functional for $\mcalS$\\
$\sd$ \> entropy per unit of mass, functional for $\mcalB$\\
$S^{\alpha}$ \> contra-variant, out-of-plane shear stress components of the Kirchhoff-Love shell\\
$\bSz$ \> surface second Piola-Kirchhoff stress tensor\\
$\bSd$ \> second Piola-Kirchhoff stress tensor\\
$\bSTd$ \> pushed forward of $\bSd$, i.e.~$\bSTd=1/\JTd\,\bFTd\,\bSd\,\bFTd^{\text{T}}$\\
$\mathcal{S}_{0}$ \>  reference configuration of the mid-surface\\
$\mathcal{S}$ \> current configuration of the mid-surface\\
$\mcalScaz$ \> current configuration of the surface through thickness\\
$\bsigd$ \> Cauchy stress tensor\\
$\bsig_{\text{KL}}$ \> Cauchy stress tensor of the Kirchhoff-Love shell\\
$\bsigz$ \> surface Cauchy stress tensor, i.e.~$\bsigz=\bFz\,\bsigz^{\sharp\lhd}\,\bFz^{\text{T}}$\\
$\bsigz^{\sharp\lhd}$ \> pull back of $\bsigz$, i.e.~$\bsigz^{\sharp\lhd}=\ds 2\bFTz^{-1}\rhoIz \pa{\psiz}{\bCez}\bFTz^{-\text{T}}$\\
$\sigab$ \> contra-variant components of $\bsigz$ ; $\alpha,\beta=1,2$\\
$t_{0}$ \> reference shell thickness\\
$t$ \> current shell thickness\\
$T$ \> absolute temperature\\
$\btId$ \> traction acting on the boundary $\partial\mcalB_{0}$ normal to $\bolds{\sV}$\\
$\btd$ \> traction acting on the boundary $\partial\mcalB$ normal to $\bnu$\\
$\btz$ \> traction acting on the boundary $\partial\mcalS$ normal to $\bnu$\\
$\tauab$ \> contra-variant components of the Kirchhoff stress tensor, i.e.~$\tauab=J\sigab$ ; $\alpha,\beta=1,2$\\
$\bwd$ \> spin of $\mcalB$\\
$\bwAd$ \> angular-velocity vector, i.e.~axial vector of $\bwd$\\
$\bXz$ \> reference position of a surface point on $\mcalS_{0}$ \\
$\bxz$ \> current position of a surface point on $\mcalS$ \\
$\bXTz$ \> intermediate position of a surface point on $\mcalSTz$ \\
$\bXcaz$ \> reference position of a surface point on $\mcalScaz_{0}$ \\
$\bxcaz$ \> current position of a surface point on $\mcalScaz$ \\
$\bXcaTz$ \> intermediate position of a surface point on $\mcalScaz_{{\text{T}}}$ \\
$\bXd$ \> reference position of a point in $\mcalB_{0}$ \\
$\bxd$ \> current position of a point in $\mcalB$ \\
$\bXTd$ \> intermediate position of a point in $\mcalBTd$ \\
\end{tabbing}
\section{Introduction}
Thermomechanical formulations for thin-walled structures have many applications. Examples from technology are sheet metal forming \citep{Bergman2004_01}, gas turbine blades \citep{Witek2016_01}, cooling systems for circuit boards \citep{Hu2017_01}, fatigue of transistors \citep{Zeanh2008_01}, batteries \citep{Guo2010_01}, solar cells \citep{Eitner2011_01}. Examples from nature are lipid membranes \citep{Sahu2017_01} and leaves \citep{Xiao2011_01}. \\
These problems can be modeled as three dimensional (3D)\footnote{By ``3D continua'' we mean 3D volumetric continua as opposed to surface continua.} or surface continua. The latter can be used to avoid the high computational costs of volumetric finite element formulations (FE), especially if they are described by Kirchhoff-Love shell kinematics. Anisotropic hyperelastic \citep{Holzapfel2000_02,Einstein2004_01,Menzel2006_01, Itskov2006_01,Gizzi2014_01,Bailly2014_01}, viscoelastic \citep{Biot1954_01,Halpin1968_01,Holzapfel2001_01,Volkov2007_01,Chen2016_01} and plastic \citep{DeBorst1990_01,Papadopoulos2001_01,Barlat2005_01,Choi2009_01,Lee2012_01,Cuomo2015_01, Li2016_01,Lee2017_01,Li2017_01} constitutive laws are well developed for volumetric continua. These material models need to be adapted to their surface counterparts in order to be used in membrane and shell formulations. Therefore the kinematics of surface and volumetric continua need to be connected to extract these sufrace material models. \citet{Steigmann2009_01,Steigmann2013_01} derives a membrane and Koiter shell formulation from 3D nonlinear elasticity. \citet{Tepole2015_01} and \citet{Roohbakhshan2016_01} utilize a projection method to extract a membrane constitutive law from the anisotropic 3D biomaterial model of \citet{Gasser2006_01}. This projection method is extended to extract a shell formulation for composite materials and biological tissues by \citet{Roohbakhshan2016_02,Roohbakhshan2017_01}.\\
There are different approaches to formulate the kinematics of deforming continua in a thermomechanical framework. The total deformation can be additively or multiplicatively decomposed into thermal and mechanical parts. The additive decomposition is widely used in the linear regime \citep{Carrera2016_01}. The linear formulation can be extended to a nonlinear one by using the multiplicative decomposition introduced by \citet{Kroner1959_01,Besseling1968_01} and \citet{Lee1969_01} for plasticity (see \ref{s:strain_measure} for a discussion on the limitations of the additive decomposition). The multiplicative decomposition of the deformation gradient is used in the current work.\\
3D thermomechanical formulations are well developed. \citet{Steinmann2002_01} presents a spatial and material framework of thermo-hyper-elastodynamics. \citet{Vujovsevic2002_01} propose a finite strain, thermoelasticity formulation based on the multiplicative decomposition of the deformation gradient. Constitutive laws for thermoelastic, elastoplastic and biomechanical materials are discussed by \citet{Lubarda2004_01}. The strong ellipticity condition should be satisfied for the stability of material models \citep[p.~258]{Marsden1994_01} and \citep{Ogden2007_01}. The stability and convexity of thermomechanical constitutive laws are considered by \citet{Silhavy1997_01} and \citet{Lubarda2008_01,Lubarda2008_02}. \citet{Miehe2011_01} propose a thermoviscoplasticity framework based on the multiplicative decomposition of elastic and plastic deformations, logarithmic strain and linear, isotropic, thermal expansion.\\
Membrane and shell formulations tend to have a much lower computational cost in comparison with volumetric formulations. \citet{Gurtin1975_01} propose a tensorial continuum theory for surface continua coupled to bulk continua. \citet{Green1979_01} propose a thermomechanical shell formulation based on Cosserat surfaces. The temperature variation through the thickness is included, and  the surface kinematics are connected to their 3D counterparts. However, the multiplicative decomposition of the deformation gradient and the effects of the thermal deformation gradient on the curvature are not considered in their work. \citet{Miehe2004_01} propose an isotropic elastoplastic solid-shell formulation and compare the additive and multiplicative decompositions of strains. \citet{Sahu2017_01} propose an irreversible thermomechanical shell formulation for lipid membranes. Since those are usually subjected to isothermal conditions, \citet{Sahu2017_01} do not consider the temperature variation, thermal expansion and heat conduction through the thickness. \citet{Steigmann2009_01,Steigmann2010_01,Steigmann2013_01} investigates the strong ellipticity and stability condition for membrane, shell and volumetric material models.\\
A shell formulation can be developed based on rotational or rotation-free formulations. \citet{Simo1989_01} propose a rotation-based, geometrically exact, shell model based on inextensible one-director and Cosserat surfaces. \citet{Simo1990_01} extend the mentioned theory to consider the thickness variation under loading. Rotation-free Kirchhoff-Love and Reissner-Mindlin shell formulations have been proposed based on subdivision surfaces \citep{Cirak2000_01,Cirak2001_01,Long2012_01}. \citet{Kiendl2009_01} propose a rotation-free, Kirchhoff-Love shell formulation using isogeometric analysis \citep{Hughes2005_01}. The formulation of \citet{Kiendl2009_01} is extended for multi patches and nonlinear materials by \citet{Kiendl2010_01,Kiendl2015_01}, respectively.  \citet{Duong2016_01} propose a rotation-free shell formulation based on the curvilinear membranes and shell formulations of \citet{Sauer2014_01} and \citet{Sauer2017_01}. \citet{Schollhammer2018_01} propose a Kirchhoff-Love shell theory based on tangential differential calculus without the introduction of a local coordinate system. See also \citet{Duong2016_01} and \citet{Schollhammer2018_01} for recent reviews of shell formulations.\\
A FE formulation is proposed by \citet{Jeffers2016_01} for heat transfer through shells. FE shell formulations for functionally graded materials are proposed by \citet{Reddy1998_01,Abbasi2000_01,Cinefra2010_01} and \citet{Kar2016_01}. \citet{Harmel2017_01} use hybrid finite element discretizations that combine isogeometric and Lagrangian interpolations for accurate and efficient thermal simulations.\\
Thermomechanical problems can be solved computationally either by monolithic or by splitting methods.
In a monolithic approach, which tends to be more robust, the mechanical and thermal parts are assembled in a single matrix and solved together. In splitting methods, the problem is divided into mechanical and thermal parts and each part is solved separately and the results are exchanged between partitions. For the partitioning, either isothermal or adiabatic splits\footnote{The adiabatic split becomes the isentropic split if there is no dissipation \citep{Holzapfel2000_01}.} can be used. The latter method can be shown to be unconditionally stable \citep{Holzapfel2000_01,Ibrahimbegovic2009_01}.\\
The condition number and accuracy of the numerical method is affected by the discretization method. The influence of the discretization on the condition number of heat conduction problems is investigated by \citet{Surana1991_01,Surana1992_01} for shells, and by \citet{Ling1994_01} for axisymmetric problems. Furthermore, the temperature and flux discontinuity should be considered in the modeling of interfaces \citep{Temizer2010_01,Temizer2014_01,Madhusudana2013_01}.\\
\citet{Sauer2018_01} present the multiplicative decomposition of the surface deformation gradient and examine its consequences on the kinematics, balance laws and constitutive relations of curved surfaces. Their formulation uses a direct surface description in curvilinear coordinates that is not necessarily connected to an underlying 3D volumetric formulation. Their examples consider isotropic material behavior. The present formulation, on the other hand, considers the derivation of surface formulations from 3D theories using a tensorial description and accounting for anisotropic material behavior.
In the current work, the kinematics of two and three dimensional continua are connected for nonisothermal materials. This connection can be used to extract membrane and shell material models from their 3D counterparts. In addition, it can be used to extend available isothermal membrane and shell material models \citep{Ghaffari2018_01,Ghaffari2018_03} to non-isothermal constitutive laws.\\
The highlights of the current work are:
\begin{itemize}
  \item The balance laws for surface and 3D continua are connected. So the extraction of the surface constitutive laws from their 3D counterparts becomes clearer.
  \item A multiplicative decomposition is used to obtain a new nonlinear thermomechanical shell formulation.
  \item Anisotropy is considered in the thermal expansion, conductivity and Helmholtz free energy.
  \item The formulation can be used to describe coupled, nonlinear thermomechanical behavior of 2D and 3D crystals, and anisotropic biological materials.
  \item The proposed thermomechanical formulation is suitable for an extension to a new thermo-electro-magneto-mechanical formulation for surface and 3D continua based on the works of \citet{Green1983_01,Chatzigeorgiou2015_01,Baghdasaryan2016_01,Dorfmann2016_01,Mehnert2017_01}.
  \item It can be used in the simulation of anisotropic thermal conductivity in printed circuit boards such as discussed in \citet{Dede2015_01,Dede2018_01}.
\end{itemize}
The remainder of this paper is organized as follows: Different descriptions of tensorial quantities and three tensorial products are introduced in Sec.~\ref{s:different_descripton_tensors}. In Sec.~\ref{s:3D_contiua}, the kinematics and equilibrium relations of a 3D anisotropic polar continua are obtained in curvilinear coordinates. Sec.~\ref{s:3D_contiua} is closed with simplified anisotropic non-polar continua. In Sec.~\ref{s:Surf_shell}, the Kirchhoff-Love shell formulation is obtained from the proposed 3D formulation. In Sec.~\ref{s:constitutive_laws}, the evolution of material symmetry groups and constitutive laws for heat transfer are discussed, and some examples of the Helmholtz free energy are presented. The paper is concluded in Sec.~\ref{s:conclusion}.

\section{Different description of tensors}\label{s:different_descripton_tensors}
In this section, some preliminary descriptions for tensorial objects in the reference and current configurations are discussed. The push forward and pull back operators are introduced for co-variant, contra-variant, and mixed tensors. Further, three tensorial products are introduced. These descriptions and tensorial products are used in the next sections for the development of a thermomechanical formulation. Tensors related to surface continua are indicated\footnote{``$\bullet$'' is a placeholder for general tensorial quantities.} by $\bullet_{\text{s}}$ while tensors related to three dimensional continua have no additional subscript. Tensors related to the intermediate configuration are written in Gothic font. The tangent vectors of three dimensional continua and of the shell mid-surface are denoted by $\bG_{i}$ and $\bA_{\alpha}$ in the reference configuration and $\bA_{\alpha}$ and $\bg_{i}$ in the current configuration. Here, Latin and Greek indexes run from 1 to 3 and 1 to 2, respectively.\\
The co-variant and contra-variant descriptions for the general vector $\bu$ in the reference configuration $\mcalB_{0}$ and the general vector $\bv$ in the current configuration $\mcalB$ are defined as\footnote{$\mcalB_{0}$ and $\mcalB$ are the reference and current configuration, see Sec.~\ref{s:3D_contiua} for details.} \citep{Marsden1994_01,Giessen1996_01,Giessen1996_02,Stumpf1997_01,Menzel2003_01,Kintzel_2006_02}
\eqb{lllll}
\ds \bud^{\sharp} \dis \ds  u^{i}\,\bG_{i}~,~~\ds \bvd^{\sharp} \dis \ds  v^{i}\,\bg_{i}~,\\[3mm]
\ds \bud^{\flat} \dis \ds  u_{i}\,\bG^{i}~,~~\ds \bvd^{\flat} \dis \ds  v_{i}\,\bg^{i}~,
\label{e:vectors}
\eqe
where
\eqb{lllll}
\ds u^{i}\is \ds \bud\cdot\bG^{i}~,~~\ds v^{i}\is \ds\bvd\cdot\bg^{i}~,\\[3mm]
\ds u_{i}\is \ds \bud\cdot\bG_{i}~,~~\ds v_{i}\is \ds \bvd\cdot\bg_{i}~,
\eqe
are the contra-variant and co-variant components of $\bud$ and $\bvd$. The co-variant, contra-variant and mixed tensors for the general second order tensors $\bU~\text{in}~\mcalB_{0}$ and $\bV~\text{in}~\mcalB$ are defined as
\eqb{lllll}
\ds \bUd^{\sharp} \dis \ds  U^{ij}\,\bG_{i}\otimes\bG_{j}~,~~\ds \bVd^{\sharp} \dis \ds  V^{ij}\,\bg_{i}\otimes\bg_{j}~,\\[3mm]
\ds \bUd^{\flat} \dis \ds  U_{ij}\,\bG^{i}\otimes\bG^{j}~,~~\ds \bVd^{\flat} \dis \ds  V_{ij}\,\bg^{i}\otimes\bg^{j}~,\\[3mm]
\ds \bUd^{\backslash} \dis \ds  U^{i}_{~j}\,\bG_{i}\otimes\bG^{j}~,~~\ds \bVd^{\backslash} \dis \ds  V^{i}_{~j}\,\bg_{i}\otimes\bg^{j}~,\\[3mm]
\ds \bUd^{/} \dis \ds  U_{i}^{~j}\,\bG^{i}\otimes\bG_{j}~,~~\ds \bVd^{/} \dis \ds  V_{i}^{j}\,\bg^{i}\otimes\bg_{j}~,
\label{e:tensors}
\eqe
where
\eqb{lllll}
\ds U^{ij} \is \bG^{i}\cdot\bUd\,\bG^{j}~,~~V^{ij} \is \bg^{i}\cdot\bVd\,\bg^{j}~,\\[3mm]
\ds U_{ij} \is \bG_{i}\cdot\bUd\,\bG_{j}~,~~V_{ij} \is \bg_{i}\cdot\bVd\,\bg_{j}~,\\[3mm]
\ds U^{i}_{~j} \is \ds \bG^{i}\cdot\bUd\,\bG_{j}~,~~ \ds V^{i}_{~j} \is \ds \bg^{i}\cdot\bVd\,\bg_{j}~,\\[3mm]
\ds \ds U_{i}^{~j}\is \ds\,\bG_{i}\cdot\bUd\,\bG^{j}~,~~\ds V_{i}^{j}\is \ds\,\bg_{i}\cdot\bVd\,\bg^{j}~,
\eqe
are the contra-variant, co-variant and mixed components of $\bUd$ and $\bVd$. Further, $\bG^{i}$ and $\bg^{i}$ are the tangent and dual vectors in the reference and current configuration (see Sec.~\ref{s:3D_contiua} for their definition). The components of surface tensors can be defined similarly, e.g.~$\bUz^{\sharp}=U^{\alpha\beta}\,\bA_{\alpha}\otimes\bA_{\beta}$ with $U^{\alpha\beta}=\bA^{\alpha}\cdot\bU\cdot\bA^{\beta}$. The vectors and tensors defined in Eq.~\eqref{e:vectors} and \eqref{e:tensors} are just different expressions of the same objects, i.e.~ $\bud=\bud^{\sharp}=\bud^{\flat}$, $\bvd=\bvd^{\sharp}=\bvd^{\flat}$, $\bUd=\bUd^{\sharp}=\bUd^{\flat}=\bUd^{\backslash}=\bUd^{/}$ and $\bVd=\bVd^{\sharp}=\bVd^{\flat}=\bVd^{\backslash}=\bVd^{/}$. The importance of these descriptions becomes clear when the push forward and pull back operators are defined.
The push forward operators for different forms of $\bUd$ are defined as
\eqb{lllll}
\ds \bUd^{\sharp\rhd} \dis \ds \bFd\,\bUd^{\sharp}\,\bFd^{\text{T}}\is U^{ij}\,\bg_{i}\otimes\bg_{j}~,\\[3mm]
\ds \bUd^{\flat\rhd} \dis \ds \bFd^{-\text{T}}\,\bUd^{\flat}\,\bFd^{-1} \is U_{ij}\,\bg^{i}\otimes\bg^{j}~,\\[3mm]
\ds \bUd^{\backslash\rhd} \dis \ds \bFd\,\bUd^{\backslash}\,\bFd^{-1} \is U^{i}_{~j}\,\bg_{i}\otimes\bg^{j}~,\\[3mm]
\ds \bUd^{/\rhd} \dis \ds \bFd^{-\text{T}}\,\bUd^{/}\,\bFd^{\text{T}} \is U_{i}^{~j}\,\bg^{i}\otimes\bg_{j}~,
\label{e:push_forward}
\eqe
and the pull back operators for different forms of $\bVd$ are defined as
\eqb{lllll}
\ds \bVd^{\sharp\lhd} \dis \ds \bFd^{-1}\,\bVd^{\sharp}\,\bFd^{-\text{T}}\is V^{ij}\,\bG_{i}\otimes\bG_{j}~,\\[3mm]
\ds \bVd^{\flat\lhd} \dis \ds \bFd^{\text{T}}\,\bVd^{\flat}\,\bFd \is V_{ij}\,\bG^{i}\otimes\bG^{j}~,\\[3mm]
\ds \bVd^{\backslash\lhd} \dis \ds \bFd^{-1}\,\bVd^{\backslash}\,\bFd \is V^{i}_{~j}\,\bG_{i}\otimes\bG^{j}~,\\[3mm]
\ds \bVd^{/\lhd} \dis \ds \bFd^{\text{T}}\,\bVd^{/}\,\bFd^{-\text{T}} \is V_{i}^{j}\,\bG^{i}\otimes\bG_{j}~,
\label{e:pull_backward}
\eqe
where $\bFd$ is the gradient deformation (see Sec.~\ref{s:3D_contiua} for details). It should be noted that the different forms of the push forward of $\bUd$ as well as the different forms for the pull back of $\bVd$, are not equal in general, i.e.~ $\bUd^{\sharp\rhd}\neq\bUd^{\flat\rhd}\neq\bUd^{\backslash\rhd}\neq\bUd^{/\rhd}$ and $\bVd^{\sharp\rhd}\neq\bVd^{\flat\rhd}\neq\bVd^{\backslash\rhd}\neq\bVd^{/\rhd}$. $\bFd$ and its transpose, inverse and transpose-inverse can be written as \citep{Basar2000_01}
\eqb{l}
\ds \bFd = \bg_{i}\otimes\bG^{i}~,~~
\bFd^{\text{T}} = \bG^{i}\otimes\bg_{i}~,~~
\bFd^{-1} = \bG_{i}\otimes\bg^{i}~,~~
\bFd^{-\text{T}} = \ds \bg^{i}\otimes\bG_{i}~.
\label{e:F_Finv_Trans}
\eqe
In Eq.~\eqref{e:pull_backward} and \eqref{e:F_Finv_Trans}, the push forward and pull back are defined using $\bFd$. They can also be defined with other two-point tensors. The push forward and pull back with the general two-point tensor $\bullet$ are indicated by superscripts $\ds \rhd(\bullet)$ and $\ds \lhd(\bullet)$. Furthermore, similar transformations can be defined for surface objects by replacing the deformation gradient $\bFd$ with the surface deformation gradient $\bFz$ (see Secs.~\ref{s:Surf_shell_Cons_laws_energy}, \ref{s:Surf_shell_Cons_laws_second_law} and \ref{s:Surf_shell_Cons_laws_weak_form}).\\
The multiplication operators\footnote{\citet{Kintzel_2006_01} and \citet{Kintzel_2006_02} use $\times$ instead of $\oplus$.} $\otimes$~, $\oplus$ and $\boxtimes$ are defined for two second order tensors $\bAd$ and $\bBd$ as
\eqb{lll}
\bAd\otimes\bBd \dis A^{ij}\,B^{kl}\,\bG_{i}\otimes\bG_{j}\otimes\bG_{k}\otimes\bG_{l}~,\\
\bAd\oplus\bBd \dis A^{ij}\,B^{kl}\,\bG_{i}\otimes\bG_{k}\otimes\bG_{l}\otimes\bG_{j} =A^{il}\,B^{jk}\,\bG_{i}\otimes\bG_{j}\otimes\bG_{k}\otimes\bG_{l}~,\\
\bAd\boxtimes\bBd \dis A^{ij}\,B^{kl}\,\bG_{i}\otimes\bG_{k}\otimes\bG_{j}\otimes\bG_{l} =A^{ik}\,B^{jl}\,\bG_{i}\otimes\bG_{j}\otimes\bG_{k}\otimes\bG_{l}~.
\eqe
The second derivative of a scaler function can then be written as
\eqb{lll}
\ds \bbCd^{\mathrm{L}} \dis \ds \frac{\partial^2{\psid}}{\partial{\bC}\oplus\partial{\bC}}= \frac{\partial^2{\psid}}{\partial{C_{ij}}\partial{C_{kl}}} \,\bG_{i}\otimes\bG_{k}\otimes\bG_{l}\otimes\bG_{j}~,
\eqe
or
\eqb{lll}
\ds \bbCd \dis \ds \frac{\partial^2{\psid}}{\partial{\bC}\otimes\partial{\bC}}= \frac{\partial^2{\psid}}{\partial{C_{ij}}\partial{C_{kl}}} \,\bG_{i}\otimes\bG_{j}\otimes\bG_{k}\otimes\bG_{l}~.
\eqe
$\bbCd$ and $\ds \bigl(\bbCd^{\mathrm{L}}\bigr)^{\text{R}}$ are connected by \citep{Kintzel_2006_01}
\eqb{lll}
\bbCd \dis \ds \bigl(\bbCd^{\mathrm{L}}\bigr)^{\text{R}}= \bigl(C^{\mathrm{L}\,ijkl}\,\bG_{i}\otimes\bG_{j}\otimes\bG_{k}\otimes\bG_{l}\bigr)^{\text{R}} =C^{\mathrm{L}\,iklj}\,\bG_{i}\otimes\bG_{j}\otimes\bG_{k}\otimes\bG_{l}~,
\eqe
where $C^{\mathrm{L}\,iklj}$ can be written as
\eqb{lll}
C^{\mathrm{L}\,iklj} \is \bG^{i}\otimes\bG^{j}:\bbCd: \bG^{k}\otimes\bG^{l}~.
\eqe
The rearrangement operators $(\bullet)^{\text{R}}$ and $(\bullet)^{\text{L}}$ cancel each other, i.e.
\eqb{lll}
\left[(\bullet)^{\text{R}}\right]^{\text{L}}=\left[(\bullet)^{\text{L}}\right]^{\text{R}}\is (\bullet)~.
\eqe
For different products of $\bAd$ and $\bBd$, $(\dots)^{\text{R}}$  can be written as \citep{Ghaffari2018_03}
\eqb{lll}
(\bAd\oplus\bBd)^{\text{R}} \is \bAd\otimes\bBd~,\\[2mm]
(\bAd\otimes\bBd)^{\text{R}} \is \bAd\boxtimes\bBd^{\text{T}}~,\\[2mm]
(\bAd\boxtimes\bBd)^{\text{R}} \is \bAd\oplus\bBd^{\text{T}}~.
\eqe
The derivative of a tensor with respect to another tensor can be written as
\eqb{lll}
\ds \frac{\partial\bAd}{\partial\bBd} \dis \ds \frac{\partial A^{ij}}{\partial B_{kl}}\,\bG_{i}\otimes\bG_{j}\otimes\bG_{k}\otimes\bG_{l}~,
\eqe
or
\eqb{lll}
\ds \frac{\partial\bAd}{\oplus\partial\bBd} \dis \ds \frac{\partial A^{ij}}{\partial B_{kl}}\,\bG_{i}\otimes\bG_{k}\otimes\bG_{l}\otimes\bG_{j}~,
\eqe
and they are connected by
\eqb{lll}
\ds \frac{\partial\bAd}{\partial\bBd} \is \ds \left(\ds \frac{\partial\bAd}{\oplus\partial\bBd}\right)^{\text{R}}~.
\eqe
\section{3D volumetric continua}\label{s:3D_contiua}
In this section, the kinematics of deforming 3D continua are discussed in a curvilinear coordinate system. Then, the conservation laws are obtained for 3D anisotropic polar continua. Finally, the relations are simplified for anisotropic non-polar thermoelastic materials.
\subsection{Curvilinear description of deforming 3D continua}
The parametric description of a 3D continua in the reference and current configuration can be written as
\eqb{lllll}
\bXd \is \bXd(\xi^{i})~,\quad \bxd \is \bxd(\xi^{i},t)~,
\eqe
where $\xi^{i}$ are the curvilinear parametric coordinates and $t$ denotes time. The tangent vectors are then given by
\eqb{lllll}
\bG_{i} \is \ds \bXd_{\!,i}~,\quad \bg_{i} \is \ds \bxd_{,i}~.
\eqe
where $\bullet_{,i}=\partial\bullet/\partial\xi^{i}$ is the parametric derivative of $\bullet$. The co-variant metric of the reference and current configuration is defined by the inner product of the tangent vectors by
\eqb{lllll}
\Guij \is \bG_{i}\cdot\bG_{j}~,\quad \guij \is \bg_{i}\cdot\bg_{j}~,
\eqe
and the contra-variant metric is connected to the co-variant metric by
\eqb{lllll}
\Gij \is [\Guij]^{-1}~,\quad \gij \is [\guij]^{-1}~,
\eqe
where $[\bullet]^{-1}$ is the inverse operator and $[\bullet]$ indicates the component matrix of $\bullet$. The dual vectors are connected to the tangent vectors by
\eqb{lllll}
\bG^{i} \is \Gij\,\bG_{j}~,\quad \bg^{i} \is \gij\,\bg_{j}~,
\eqe
and the inner product of the tangent and dual vectors are
\eqb{lllll}
\delta_{i}^{j} \is \bG_{i}\cdot\bG^{j}~,\quad \delta_{i}^{j} \is \bg_{i}\cdot\bg^{j}~,
\eqe
where $\delta_{i}^{j}$ is Kronecker delta.
The 3D identity tensor $\bolds{1}$ can be written as
\eqb{lllllllll}
\bolds{1} \is \Guij\,\bG^{i}\otimes\bG^{j} \is \Gij\,\bG_{i}\otimes\bG_{j} \is \bG^{i}\otimes\bG_{i}\is \bG_{i}\otimes\bG^{i}\\
\eqe
and
\eqb{lllllllll}
\bolds{1} \is \guij\,\bg^{i}\otimes\bg^{j} \is \gij\,\bg_{i}\otimes\bg_{j} \is \bg^{i}\otimes\bg_{i}\is \bg_{i}\otimes\bg^{i}~,
\eqe
such that $\bUd=\bolds{1}\cdot\bUd=\bUd\cdot\bolds{1}$ and $\bVd=\bolds{1}\cdot\bVd=\bVd\cdot\bolds{1}$. The gradient operator for a general scalar, vector (first order tensor) and second order tensor in the reference and current configuration are
\eqb{lllll}
\Gradd\, \Phi \dis \Phi_{,i}\,\bG^{i}~,\\[3mm]
\Gradd\, \bud \dis \bud_{,j}\otimes\bG^{j} \is u^{i}_{\|j}\,\bG_{i}\otimes\bG^{j}~,\\[3mm]
\Gradd\, \bUd \dis \bUd_{\!,k}\otimes\bG^{k} \is U^{ij}_{~\|k}\,\bG_{i}\otimes\bG_{j}\otimes\bG^{k}~\\
\eqe
and
\eqb{lllll}
\gradd\, \Phi \dis \Phi_{,i}\,\bg^{i}~,\\[3mm]
\gradd\, \bvd \dis \bvd_{,j}\otimes\bg^{j} \is v^{i}_{|j}\,\bg_{i}\otimes\bg^{j}~,\\[3mm]
\gradd\, \bVd \dis \bVd_{\!,k}\otimes\bg^{k} \is V^{ij}_{~~|k}\,\bg_{i}\otimes\bg_{j}\otimes\bg^{k}~,\\
\eqe
where ``$\|$'' and ``$|$'' are the 3D co-variant derivatives in the reference and current configuration, defined by \citep{Itskov2015_01}
\eqb{lll}
u^{i}_{\|j} \dis u^{i}_{,j}+u^{k}\,\GammaId{k}{j}{i}~,\\[3mm]
U^{ij}_{~\,\|k} \dis U^{ij}_{~\,,k}+U^{lj}\,\GammaId{l}{k}{i}+U^{il}\,\GammaId{l}{k}{j}~,\\[3mm]
v^{i}_{|j} \dis v^{i}_{,j}+v^{k}\,\Gammad{k}{j}{i}~,\\[3mm]
V^{ij}_{~~|k} \dis V^{ij}_{~~,k}+V^{lj}\,\Gammad{l}{k}{i}+V^{il}\,\Gammad{l}{k}{j}~.\\[3mm]
\eqe
Here, $\Gammad{i}{j}{k}$ and $\GammaId{i}{j}{k}$ are the 3D Christoffel symbols of the second kind in the reference and current configuration, defined as
\eqb{lll}
\GammaId{i}{j}{k} \dis \bG_{i,j}\cdot\bG^{k}~
\eqe
and
\eqb{lll}
\Gammad{i}{j}{k} \dis \bg_{i,j}\cdot\bg^{k}~.
\eqe
The 3D co-variant derivative of the tangent and dual vector, and the metric in the reference and current configuration are zero (Ricci's theorem \citep[p.~54]{Itskov2015_01}), i.e.
\eqb{llllll}
\bG_{i\|j} \is \bG_{i,j}-\bG_{l}\,\GammaId{i}{j}{l}=\bolds{0}~,\quad & \bG^{i}_{\|j} \is \bG^{i}_{,j}+\bG^{l}\,\GammaId{l}{j}{i}=\bolds{0}~,\\[3mm]
G_{ij\|k} \is 0~,\quad & G^{ij}_{~\|k}\is 0
\label{e:3D_Ricci_ref}
\eqe
and
\eqb{llllll}
\bg_{i|j} \is \bg_{i,j}-\bg_{l}\,\Gammad{i}{j}{l}=\bolds{0}~,\quad & \bg^{i}_{|j} \is \bg^{i}_{,j}+\bg^{l}\,\Gammad{l}{j}{i}=\bolds{0}~,\\[3mm]
g_{ij|k} \is 0~,\quad & g^{ij}_{~\,|k}\is 0~.
\label{e:3D_Ricci_cur}
\eqe
Using Ricci's theorem, the 3D co-variant derivative can be represented in other arrangements of indices, e.g.
\eqb{lll}
\ds U_{ij\|}^{~~~k} \is \ds G_{il}\,G_{jm}\,G^{kn}\,U^{lm}_{~~\|n}~,\\[3mm]
\ds V_{ij|}^{~~k} \is \ds g_{il}\,g_{jm}\,g^{kn}\,V^{lm}_{~~~|n}~,
\eqe
where $U_{ij\|}^{~~~k}$ and $V_{ij|}^{~~k}$ are the contra-variant derivatives. These transformations can be used to define other expressions of the gradient operator. The divergence of a vector and a second order tensor are defined by
\eqb{llllll}
\Divd\, \bu \dis \bu_{,j}\cdot\bG^{j} = u^{i}_{\|i}~,\quad & \divd\, \bv \dis \bv_{,j}\cdot\bg^{j} = v^{i}_{|i}~,\\[3mm]
\Divd\, \bU \dis \bU_{\!,j}\cdot\bG^{j} = U^{ij}_{~\|j}\,\bG_{i}~,\quad & \divd\, \bV \dis \bV_{\!,j}\cdot\bg^{j} = V^{ij}_{~~|j}\,\bg_{i}~.\\[3mm]
\eqe
With this, the divergence theorem can be used to transfer a domain integral to a boundary integral as
\eqb{lll}
\ds \int_{\partial\sB_{0}}{\bud\cdot\bolds{\sV}\,\dif S} \is \ds\int_{\sB_{0}}{\Divd\,\bud\,\dif V}~,\\[3mm]
\ds \int_{\partial\sB_{0}}{\bUd\,\bolds{\sV}\,\dif S} \is \ds\int_{\sB_{0}}{\Divd\,\bUd\,\dif V}~\\[3mm]
\eqe
and
\eqb{lll}
\ds \int_{\partial\sB}{\bvd\cdot\bnu\,\dif s} \is \ds\int_{\sB}{\divd\,\bvd\,\dif v}~,\\[3mm]
\ds \int_{\partial\sB}{\bVd\,\bnu\,\dif s} \is \ds\int_{\sB}{\divd\,\bVd\,\dif v}~,\\[3mm]
\eqe
where $\sB_{0}$ and $\sB$ are the reference and current configuration of 3D continua and $\partial\sB_{0}$ and $\partial\sB$ are their corresponding boundaries. $\bolds{\sV}$ and $\bnu$ are the normal unit vectors on $\partial\sB_{0}$ and $\partial\sB$.
\subsection{Kinematics of deformation}
\begin{figure}[h]
        \centering
      \includegraphics[height=70mm]{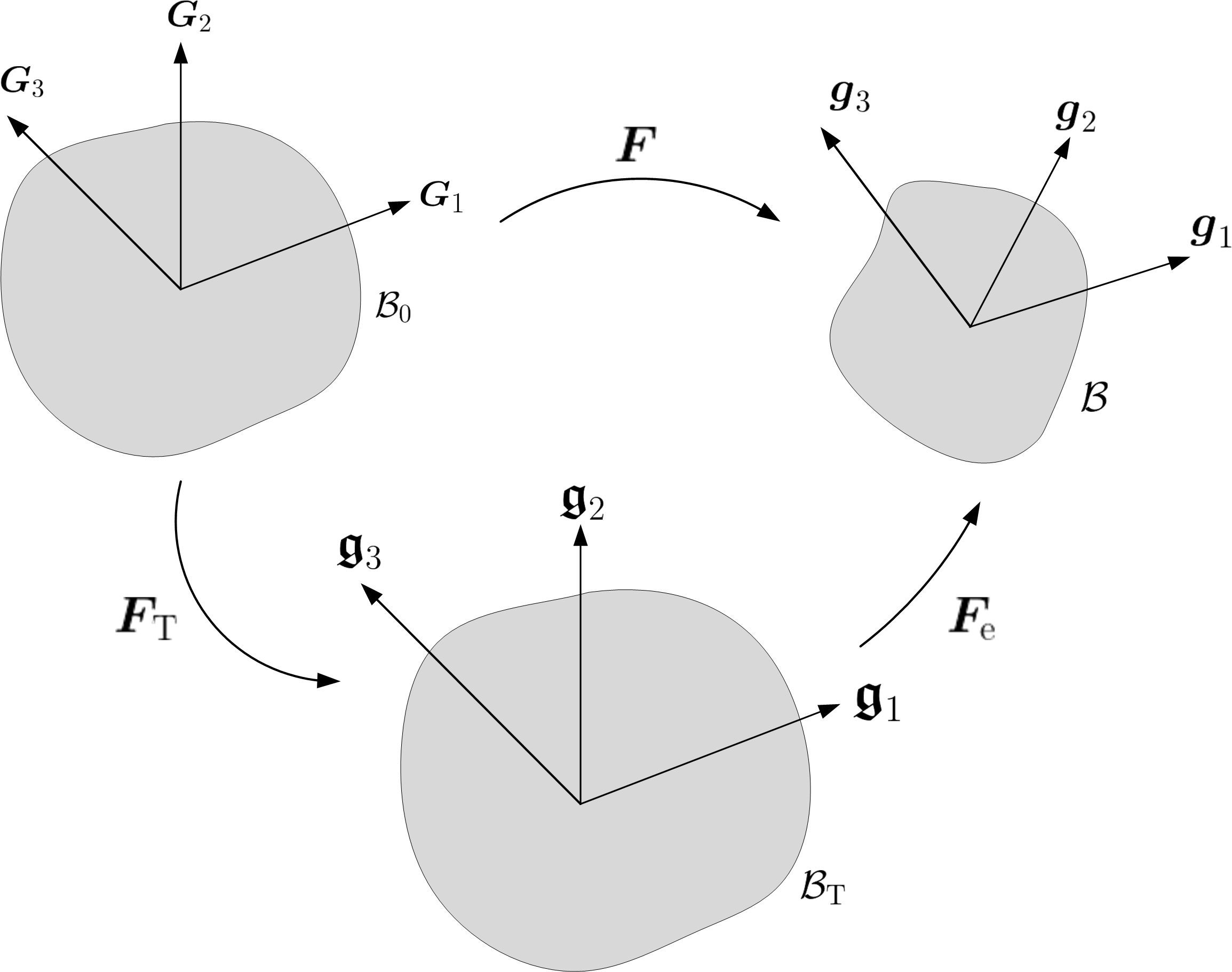}
      \caption{Multiplicative deformation decomposition of 3D volumetric continua. $\mcalB_{0}$, $\mcalBTd$ and $\mcalB$ are the reference, intermediate and current configuration and their corresponding tangent vectors are $\bG_{i}$, $\bgTdui{i}$ and $\bg_{i}$. $\bFd$, $\bFTd$ and $\bFed$ are the total, thermal and elastic deformation gradient. \label{f:3D_thermal_kin}}
\end{figure}

Next, the kinematics of deformation is discussed. An intermediate configuration $\sB_{\text{T}}$ is introduced\footnote{The deformation of the intermediate configuration can be incompatible \citep{bonet_2008_01}.}, and the reference, intermediate and current configuration $\sB_{0}$, $\sB_{\text{T}}$ and $\sB$ are connected (see Fig.~\ref{f:3D_thermal_kin}).
\noindent
An additive decomposition of the strain should not be used in finite deformations, unless the logarithmic strain is used (see \ref{s:strain_measure}). The limitation of the additive decomposition of strains is resolved by using a multiplicative decomposition of the deformation gradient. $\sB_{0}$ is mapped to $\sB$ by $\bxd=\bvarphid(\bXd)$ and the deformation gradient can be written as (see also Eq.~\eqref{e:F_Finv_Trans})
\eqb{lll}
\bFd \is \ds \pa{\bxd}{\bXd}= \pa{\bxd}{\xi^{i}}\otimes\pa{\xi^{i}}{\bXd}=\bg_{i}\otimes~\bG^{i}~.
\eqe
It pushes forward the differential line $\dif\bXd$ to $\dif\bxd$ as
\eqb{lll}
\dif\bxd \is \ds \bFd\,\dif\bXd~,
\label{3D_kin_def_grad}
\eqe
and hence pushes forward $\bG_{i}$ to $\bg_{i}$ as
\eqb{lll}
\bg_{i} \is \ds \bFd\,\bG_{i}~.
\eqe
The deformation gradient $\bFd$ can be decomposed into the elastic and thermal parts $\bFed$ and $\bFTd$ as
\eqb{lll}
\ds \bFd \is \ds \bFed\,\bFTd~.
\eqe
$\bFTd$ pushes forward $\dif\bX$ to $\dif\bXTd$ as
\eqb{lll}
\dif\bXTd \is \bFTd\,\dif\bXd~,
\eqe
and hence pushes forward $\bG_{i}$ to $\bgTdui{i}$ as
\eqb{lll}
\bgTdui{i} \is \bFTd\,\bG_{i}~,
\eqe
where $\dif\bXTd$ and $\bgTdui{i}$ are the differential line and the tangent vectors of $\sB_{\text{T}}$, respectively, and $\bgTdui{i}$ are defined as
\eqb{lll}
\bgTdui{i} \dis \bXTd_{,i}~.
\eqe
$\bFed$ pushes forward $\dif\bXTd$ to $\dif\bx$ as
\eqb{lll}
\dif\bxd \is \bFed\,\dif\bXTd~,
\eqe
and hence pushes forward $\bgTdui{i}$ to $\bg_{i}$ as
\eqb{lll}
\bg_{i} \is \bFed\,\bgTdui{i}~.
\eqe
The right Cauchy-Green deformation tensor $\bCd$ and its thermal, $\bCTd$, and elastic, $\bCed$, parts can be defined as
\eqb{lll}
\bCd \dis \bFd^{\text{T}}\,\bFd=\bFTd^{\text{T}}\,\bCed\,\bFTd~,\\[3mm]
\bCTd \dis \bFTd^{\text{T}}\,\bFTd~,\\[3mm]
\bCed \dis \bFed^{\text{T}}\,\bFed~.
\label{e:3D_right_cauchy_diff}
\eqe
The spatial velocity gradient $\bld$ can be written as
\eqb{lll}
\ds \bld \is \ds \pa{\bvd}{\bxd}=\pa{\bvd}{\bXd}\,\pa{\bXd}{\bxd} =\dot{\bFd}\,\bFd^{-1}~,
\eqe
where $\dot{\bullet}$ denotes the material time derivative of $\bullet$, i.e.
\eqb{lll}
\ds \dot{\bullet} \is \ds \frac{\text{D} \bullet}{\text{D} t}=\left.\pa{\bullet}{t}\right|_{\bXd=\text{fixed}}~,
\eqe
and the velocity $\bvd$ is the material time derivative of $\bxd$, i.e.~$\bvd=\dot{\bxd}$.
$\bld$ can be decomposed into a symmetric (the rate of deformation tensor $\bdd$) and skew symmetric (the spin tensor $\bwd$) part as
\eqb{lll}
\ds \bld \is \ds \bdd+\bwd~.
\eqe
It can also be written in terms of the elastic and thermal velocity gradients $\bld_{\text{e}}$ and $\bld_{\text{T}}$ as
\eqb{lll}
\bld \is \bled+
\bFed\,\blTd\,\bFed^{-1}~,
\eqe
with
\eqb{lll}
\bled \is \dot{\bFd}_{\text{e}}\,\bFed^{-1}~,
\eqe
\eqb{lll}
\blTd \is \dot{\bFd}_{\text{T}}\,\bFTd^{-1}~.
\eqe
The angular-velocity vector\footnote{Axial vector of $\bwd$.} $\bwAd$ and spin $\bwd$ can be connected by \citep{Malvern1969_01,Basar2000_01,bonet_2008_01}
\eqb{lll}
\bwAd \is \ds \frac{1}{2}\curld\,\bvd=\frac{1}{2}\bsEd:\left(\gradd\,\bvd\right)^{\text{T}} =\frac{1}{2}\bsEd:\bld^{\text{T}}= \frac{1}{2}\bsEd:\bwd^{\text{T}}~,
\label{e:ang_spin_to_vel}
\eqe
\eqb{lll}
\bwd^{\text{T}} \is \bwAd\cdot\bsEd~,
\label{e:ang_vel_to_spin}
\eqe
where $\bsEd$ is the 3D permutation (or alternating or Levi-Civita) tensor \citep{Basar2000_01,Holzapfel2000_01} and it can be written in Cartesian coordinates as
\eqb{lll}
\bsEd \is \sE^{ijk}\,\be_{i}\otimes\be_{j}\otimes\be_{k} = \sE_{ijk}\,\be^{i}\otimes\be^{j}\otimes\be^{k}~,
\eqe
with
\eqb{lll}
\sE_{ijk}=\sE^{ijk} \is \left\{
\begin{array}{l l}
                 1~~~&;~~~\mathrm{if}~~~i,~j,~k~\text{is even permutation of 1, 2, 3}~,\\
                 -1~~~&;~~~\mathrm{if}~~~i,~j,~k~\text{is odd permutation of 1, 2, 3}~,\\
                 0~~~&;~~~\mathrm{else}~,
\end{array}\right .
\eqe
where $\be_{i}=\be^{i}$ are the Cartesian basis vectors. The following relation can be used to express the cross product of two vectors $\bud$ and $\bvd$\footnote{See \citet{Basar2000_01} and \citet{bonet_2008_01}.}
\eqb{lll}
\bud\times\bvd \is \bsEd:(\bud\otimes\bvd)~,
\eqe
where $\times$ denotes the cross product. The $\curld$ of $\bvd$ can then be written as \citep{Green1968_01,Basar2000_01}
\eqb{lll}
\curld\,\bvd \dis  v_{i|j}\,\bg^{j}\times\bg^{i}=\bsEd:\left(\gradd\,\bvd\right)^{\text{T}}~.
\eqe
A volume element of the reference configuration $\dif V$ can be connected to its counterpart in the current configuration $\dif v$ by
\eqb{lll}
\ds \Jd \dis \ds \frac{\dif v}{\dif V} = \detd \bFd~,
\label{e:3Dvolume_conv}
\eqe
where the determinant is defined as
\eqb{lll}
\ds \detd\bullet \dis \ds \frac{[(\bullet\cdot\by_{1})\times(\bullet\cdot\by_{2})]\cdot(\bullet\cdot\by_{3})} {[\by_{1}\times\by_{2}]\cdot\by_{3}}~,
\eqe
for three non-coplanar vectors $\by_{i}$. The surface element in the reference and current configuration are connected by Nanson's formula \citep{bonet_2008_01}
\eqb{lll}
\bnu\,\dif s \is \Jd\,\bFd^{-\text{T}}\cdot\bolds{\sV}\,\dif S~.
\label{e:3Darea_conv}
\eqe
\subsection{Conservation laws}
Here, the global and local form of the conservation of mass, linear momentum, angular momentum, and the first and second law of thermodynamics are stated.
\subsubsection{Mass balance}
Assuming mass conservation and using Eq.~\eqref{e:3Dvolume_conv}, the mass balance in the reference and current configuration can be connected by
\eqb{lll}
\ds \int_{\sB}{\dif m} \is \ds \int_{\sB}{\rhod\,\dif v} =\int_{\sB_{0}}{\rhod\,\Jd\,\dif V}= \int_{\sB_{0}}{\rhoId\,\dif V}~.
\eqe
where $\rhoId$ and $\rhod$ are the mass density of the reference and current configuration. Hence, the local relation for the conservation of mass is
\eqb{lll}
\rhoId \is \Jd\,\rhod~.
\label{e:3D_cons_mass_loc}
\eqe
\subsubsection{Linear momentum balance}
The linear momentum balance can be written as
\eqb{lll}
\ds \frac{\text{D}}{\text{D}t}{\int_{\sB}{\rhod\,\bvd\,\dif v}} \is \ds \int_{\partial_{\btd}\sB}{\btd\,\dif s}+ \int_{\sB}{\rhod\,\bff\,\dif v}~,
\label{e:3D_cons_mom_int}
\eqe
where $\bff$ and $\bvd$ are the body force per unit mass and the velocity, respectively. $\btd$ and $\partial_{\btd}\sB$ are the traction (force per unit area) and its corresponding boundary.
Using the divergence theorem and Cauchy's formula
\eqb{lll}
\btd=\bsigd^{\text{T}}\,\bnu~,
\label{e:cauchy_theorem}
\eqe
Eq.~\eqref{e:3D_cons_mom_int} can be written in its local current form as
\eqb{lll}
\ds \rhod\,\dot{\bv} \is \ds \divd\,\bsigd^{\text{T}} + \rhod\,\bff~
\label{e:con_mon_loc_cur}
\eqe
and in its local reference form as
\eqb{lll}
\ds \rhoId\,\dot{\bv} \is \ds \Divd\,\bPd^{\text{T}} + \rhoId\,\bff~.
\eqe
Using Eq.~\eqref{e:3Darea_conv}, the 3D first Piola Kirchhoff stress tensor $\bPd$ can be obtained from
\eqb{lll}
\ds \int_{\partial_{\bt}\sB}{\btd\,\dif s} \is \ds \int_{\partial_{\btd_{0}}\sB}{\btd_{0}\,\dif S}~,
\eqe
as
\eqb{lll}
\ds \bPd \is \ds \Jd\bFd^{-1}\,\bsigd~,
\eqe
where $\btd_{0}=\bPd^{\text{T}}\,\bolds{\sV}$. The 3D second Piola Kirchhoff stress is defined as
\eqb{lll}
\ds \bSd \dis \ds \bPd\,\bFd^{-\text{T}}=\Jd\,\bFd^{-1}\,\bsigd\,\bFd^{-\text{T}}~.
\eqe
\subsubsection{Angular momentum balance}
If there is no spin angular momentum (intrinsic angular momentum) in continua, the angular momentum balance can be written as \citep[p.~218]{Malvern1969_01}
\eqb{lll}
\ds \frac{\text{D}}{\text{D}t}{\int_{\sB}{\bxd\times\rhod\,\bvd\,\dif v}} \is \ds \int_{\sB}{\rhod\,\left[\bxd\times\,\bff\,+\bc\right]\,\dif v} + \int_{\partial_{\btd}\sB}{\left[\bxd\times\btd+\bmd\right]\,\dif s}~,
\label{e:3D_cons_ang_mom_int}
\eqe
where $\bc$ is the body force couple (per unit mass) and $\bmd$ is the traction couple (per unit current area). Eq.~\eqref{e:3D_cons_ang_mom_int} can be simplified using $\bv\times\bv=\bzero$,\footnote{$\bar{\bullet}$ is used to distinguish between $\bmurz$ and $\bmuz$ ($\bmurz^{\text{T}}=\bn\times\bmuz^{\text{T}}$, see \ref{s:Surf_shell_Cons_laws_angular}).}
\eqb{lll}
\bmd=\bmurd^{\text{T}}\,\bnu~,
\label{e:moment_theorem}
\eqe
and\footnote{$\bvd$ and $\bVd$ are two general first and second order tensors.} \citep[p.~54]{Holzapfel2000_01}
\eqb{lll}
\divd(\bvd\times\bVd) \is \bvd\times\divd\,\bVd+\bsEd:\left[\text{grad}\,(\bvd)\,\bVd^{\text{T}}\right]~,
\label{e:div_vec_cross_ten}
\eqe
as
\eqb{lll}
\ds {\int_{\sB}{\bxd\times\rhod\,\dot{\bvd}\,\dif v}} \is \ds \int_{\sB}{\rhod\,\left[\bxd\times\,\bff\,+\bc\right]\,\dif v}+ \int_{\sB}{\left[\bxd\times\divd\,\bsigd^{\text{T}}+ \bsEd:\bsigd+\divd\,\bmurd^{\mathrm{T}}\right]\,\dif v}~.
\eqe
Here, $\bmurd$ is the moment tensor that is analogous to the Cauchy stress tensor. Using Eq.~\eqref{e:con_mon_loc_cur}, this relation can be simplified and written in its local form as
\eqb{lll}
\divd\,\bmurd^{\text{T}}+\rhod\,\bc+\bsEd:\bsigd\is \bzero~.
\label{e:con_ang_mon_loc_cur}
\eqe
See \citep[p.~220]{Malvern1969_01} for the componentwise expression of Eq.~\eqref{e:con_ang_mon_loc_cur} in Cartesian coordinates. If $\bmurd$ and $\bc$ are zero, $\bsigd^\text{T}= \bsigd$ follows.
\subsubsection{The first law of thermodynamics}
The energy balance in its integral form for 3D continua can be written as
\eqb{lll}
\ds \frac{\text{D}}{\text{D}t}\left(\sK +\sU\right) \is \ds \sP_{\text{ext}}+\sQ~,
\label{e:3D_cons_first_law_int}
\eqe
where $\sK$ , $\sU$ , $\sP_{\text{ext}}$ and $\sQ$ are the kinetic energy, internal energy, rate of external mechanical power and rate of heat input that can be written as \citep[p.~220]{Malvern1969_01}
\eqb{lll}
\ds \sK \is \ds \int_{\sB}{\frac{1}{2}\rhod\,\bvd\cdot\bvd\,\dif v}~,
\eqe
\eqb{lll}
\ds \sU \is \ds \int_{\sB}{\rhod\,\ud\,\dif v}~,
\eqe
\eqb{lll}
\sP_{\text{ext}} \is \ds \int_{\sB}{\rhod\,\left[\,\bff\cdot\bvd+\bc\cdot\bwAd\right]\,\dif v} + \int_{\partial_{\btd\bmd}\sB}{\left[\btd\cdot\bvd+\bmd\cdot\bwAd\right]\,\dif s}~,
\eqe
\eqb{lll}
\sQ \is \ds \int_{\sB}{\rhod\,r\,\dif v}-\int_{\partial_{\bqd}\sB}{\bqd\cdot\bnu\,\dif s}~,
\eqe
where $\bwAd$ is defined in Eq.~\eqref{e:ang_spin_to_vel}. Here, $\bqd$, $r$ and $\ud$ are the heat flux vector per unit area, heat source and internal energy per unit mass, respectively. Using these relations, Eq.~\eqref{e:3D_cons_first_law_int} can be written as
\eqb{lll}
\ds \int_{\sB}{\rhod\,\dot{\bvd}\cdot\bvd\,\dif v}+\int_{\sB}{\rhod\,\dot{\ud}\,\dif v} \is
\ds \int_{\sB}{\rhod\,\left[\,\bff\cdot\bvd+\bc\cdot\bwAd\right]\,\dif v} + \int_{\partial_{\btd\bmd}\sB}{\left[\btd\cdot\bvd+\bmd\cdot\bwAd\right]\,\dif s}\\[4mm] \plus \ds \int_{\sB}{\rhod\,r\,\dif v}-\int_{\partial_{\bqd}\sB}{\bqd\cdot\bnu\,\dif s}~.
\label{e:con_eng_int_cur_form1}
\eqe
Using the divergence theorem, and Eqs.~\eqref{e:cauchy_theorem}, \eqref{e:con_mon_loc_cur} and \eqref{e:moment_theorem}, $\sP_{\text{ext}}$ can be written as
\eqb{lll}
\sP_{\text{ext}}\is\ds \int_{\sB}{\left[\left(\divd\,\bsigd^{\text{T}}+\rhod\,\bff\right)\cdot\bvd+ \bsigd^{\text{T}}:\bld + \left(\divd\,\bmurd^{\text{T}}+\rhod\,\bc\right)\cdot\bwAd+ \bmurd^{\text{T}}:\gradd\,\bwAd\right]\,\dif v}~.
\eqe
Using this relation and Eq.~\eqref{e:con_ang_mon_loc_cur}, Eq.~\eqref{e:con_eng_int_cur_form1} can be rewritten as
\eqb{lll}
\ds \int_{\sB}{\rhod\,\dot{\ud}\,\dif v} \is
\ds \int_{\sB}{\left[\bsigd^{\text{T}}:\bld - \bwAd\cdot\bsEd:\bsigd+\bmurd^{\text{T}}:\gradd\,\bwAd +\rhod\,r-\divd\,\bqd\right]\,\dif v}~.
\label{e:con_eng_int_cur_form1_1}
\eqe
Next, $\bsigd$ is decomposed into symmetric and skew symmetric parts as
\eqb{lll}
\bsigd\is \bsigd_{\text{sym}}+\bsigd_{\text{skew}}~.
\eqe
Substituting this relation and Eq.~\eqref{e:ang_vel_to_spin} into $\bsigd^{\text{T}}:\bld$ results in
\eqb{lll}
\bsigd^{\text{T}}:\bld \is \bsigd_{\text{sym}}^{\text{T}}:\bdd+\bsigd_{\text{skew}}^{\text{T}}:\bwd\\[4mm]
\is  \bsigd_{\text{sym}}^{\text{T}}:\bdd+\bwAd\cdot\bsEd:\bsigd~.
\eqe
Using this relation, Eq.~\eqref{e:con_eng_int_cur_form1_1} can be simplified as
\eqb{lll}
\ds \int_{\sB}{\left(\rhod\,\dot{\ud}-\bsigd^{\text{T}}_{\text{sym}}:\bdd -\bmurd^{\text{T}}:\gradd\,\bwAd-\rhod\,r +\divd\,\bqd\right)\,\dif v} \is 0~,
\label{e:con_eng_int_cur_form1_2}
\eqe
so that its local form becomes
\eqb{lll}
\ds \rhod\,\dot{\ud}-\bsigd^{\text{T}}_{\text{sym}}:\bdd- \bmurd^{\text{T}}:\gradd\,\bwAd-\rhod\,r+\divd\,\bqd \is 0~.
\label{e:con_eng_loc_cur_form1}
\eqe
An alternative form for Eq.~\eqref{e:con_eng_int_cur_form1_2} can be obtained by using the divergence theorem, Eq.~\eqref{e:cauchy_theorem}, \eqref{e:con_mon_loc_cur}, \eqref{e:moment_theorem}, \eqref{e:con_ang_mon_loc_cur} and \eqref{e:con_eng_int_cur_form1}. This gives
\eqb{lll}
\ds \int_{\sB}{\rhod\,\dot{\ud}\,\dif v} \is
\ds \int_{\sB}{\left[\bsigd^{\text{T}}:\bld+\divd\,\left(\bwAd\cdot\bmurd^{\text{T}}\right) +\rhod\,\left(r+\bc\cdot\bwAd\right)-\divd\,\bqd\right]\,\dif v}~,
\eqe
or in local form
\eqb{lll}
\rhod\,\dot{\ud} -\bsigd^{\text{T}}:\bld-\divd\,\left(\bwAd\cdot\bmurd^{\text{T}}\right) -\rhod\,\left(r+\bc\cdot\bwAd\right)+\divd\,\bqd \is 0~.
\label{e:con_eng_loc_cur_form2}
\eqe
The first law in the reference configuration can be written as
\eqb{lll}
\rhoId\,\dot{\ud} -\bPd^{\text{T}}:\dot{\bFd}-\Divd\,\left(\bwAd\cdot\bmurd^{\text{T}}_{0}\right) -\rhoId\,\left(r+\bc\cdot\bwAd\right)+\Divd\,\bQd \is 0~.
\label{e:con_eng_loc_ref_form2}
\eqe
From
\eqb{lll}
\ds \int_{\partial_{\bqd}\sB_{0}}{\bQd\cdot\bolds{\sV}\,\dif S} \is \ds \int_{\partial_{\bqd}\sB}{\bqd\cdot\bnu\,\dif s }= \int_{\partial_{\bqd}\sB_{0}}{\bqd\,\Jd\bFd^{-\text{T}}\,\bolds{\sV}\,\dif S}~,
\eqe
the relationship between $\bQd$ and $\bqd$ is found as
\eqb{lll}
\bQd \is \bqd\,\Jd\,\bFd^{-\text{T}}~,
\eqe
and from
\eqb{lll}
\ds \int_{\partial_{\bmd}\sB_{0}}{\bwAd\cdot\bmurd^{\text{T}}_{0}\,\bolds{\sV}\,\dif S} \is \ds \int_{\partial_{\bmd}\sB}{\bwAd\cdot\bmurd^{\text{T}}\,\bnu\,\dif s }= \int_{\partial_{\bmd}\sB_{0}}{\bwAd\cdot\bmurd^{\text{T}}\,\Jd\bFd^{-\text{T}}\,\bolds{\sV}\,\dif S}~,
\eqe
the relationship between $\bmurd_{0}$ and $\bmurd$ is found as
\eqb{lll}
\ds \bmurd_{0} \is \ds \Jd\,\bFd^{-1}\bmurd~.
\eqe
\subsubsection{The second law of thermodynamics}
Here, the second law of thermodynamics in form of the Clausius-Duhem inequality is obtained. Then, the Clausius-Planck inequality for the case of non-negative internal dissipation is derived from the Clausius-Duhem inequality. It should be emphasized that ``the truth of the Clausius-Planck inequality should not be assumed. Rather, it can sometimes be proved to hold as a theorem'' \citep{Truesdell1967_01}. The second law of thermodynamics can be written in the integral form as
\eqb{lll}
\ds \frac{\text{D}}{\text{D}t}\int_{\sB}{\rhod\,\sd\,\dif v}\geq
\ds \int_{\sB}{\rhod\,\frac{r}{T}\,\dif v} - \int_{\partial_{q}\sB}{\frac{\bqd\cdot\bnu}{T}\,\dif s}
\eqe
and in the local form as
\eqb{lll}
\ds \rhod\,\dot{\sd}-\rhod\,\frac{r}{T}+\divd\left(\frac{\bqd}{T}\right)\geq 0~,
\eqe
where $\sd$ is the entropy (per mass) and $T$ is the absolute temperature in Kelvin. Using the latter relation and Eq.~\eqref{e:con_eng_loc_cur_form1}, the Clausius-Duhem inequality can be obtained as
\eqb{lll}
\ds \rhod\,T\,\dot{\sd}-\rhod\,\dot{\ud}+\bsigd_{\text{sym}}^{\text{T}}:\bdd+ \bmurd^{\text{T}}:\gradd\,\bwAd-\frac{1}{T}\bqd\cdot\gradd(T) \geq 0~.
\eqe
Alternatively, by using Eq.~\eqref{e:con_eng_loc_cur_form2}, it can be written in the current configuration as
\eqb{lll}
\ds \rhod\,T\,\dot{\sd}-\rhod\,\dot{\ud}+\bsigd^{\text{T}}:\bld+ \divd\,\left(\bwAd\cdot\bmurd^{\text{T}}\right)+\rhod\,\bc\cdot\bwAd -\frac{1}{T}\bqd\cdot\gradd(T) \geq 0~,
\label{e:con_second_law_loc_cur_form}
\eqe
and in the reference configuration as
\eqb{lll}
\ds \rhoId\,T\,\dot{\sd}-\rhoId\,\dot{\ud}+\bPd^{\text{T}}:\dot{\bFd}+ \Divd\,\left(\bwAd\cdot\bmurId^{\text{T}}\right)+ \rhoId\,\bc\cdot\bwAd-\frac{1}{T}\bQd\cdot\Gradd(T) \geq 0~,
\eqe
This relation includes the rate of the local and conductive entropy production $\gamma_{\text{loc}}$ and $\gamma_{\text{con}}$ that can be written as
\eqb{lll}
\ds \gamma_{\text{loc}} \is \ds \rhod\,\dot{\sd}-\rhod\,\frac{r}{T}+\frac{1}{T}\divd(\bqd)= \rhod\,\dot{\sd}-\frac{\rhod\,\dot{\ud}}{T}+\frac{1}{T}\left[\bsigd^{\text{T}}:\bld+ \divd\,\left(\bwAd\cdot\bmurd^{\text{T}}\right)+\rhod\,\bc\cdot\bwAd\right]~,\\[3mm]
\ds \gamma_{\text{con}} \is \ds -\frac{1}{T^2}\bqd\cdot\gradd\,T~.
\label{e:con_second_law_loc_cur_form2}
\eqe
The Clausius-Planck inequality requires that $\gamma_{\text{loc}}\geq 0$ and $\gamma_{\text{con}} \geq 0$ \citep{Truesdell2004_01,Mahnken2008_01}.
\subsection{3D non-polar continua}
If the body force couple $\bc$ and traction moment $\bmd$ are assumed to be zero, the local conservation laws can be summarized as
\eqb{c}
\ds \rhoId = \Jd\,\rhod~.\\
\ds \rhoId\,\dot{\bvd} =  \Divd\,\bPd^{\text{T}} + \rhoId\,\bff~, \\
\ds \bsigd= \bsigd^{\text{T}}~,\\
\ds \rhoId\,\dot{\ud}=\frac{1}{2}\bSd:\dot{\bCd} +\rhoId\,r-\Divd\,\bQd~,\\
\ds \rhoId\,T\,\dot{\sd}-\rhoId\,\dot{\ud}+\frac{1}{2}\bSd:\dot{\bCd} -\frac{1}{T}\bQd\cdot\Gradd(T) \geq 0~.
\label{e:con_laws_3D_non_polar}
\eqe
It is more convenient to write the second law of thermodynamics in terms of the Helmholtz free energy per mass
\eqb{lll}
\psid \is \ud-T\,\sd~,
\eqe
which has the time derivative
\eqb{lll}
\dot{\psid} \is \dot{\ud}-\dot{T}\,\sd-T\,\dot{\sd}~.
\eqe
Substitution of this relation into Eq.~\eqref{e:con_laws_3D_non_polar} results in
\eqb{lll}
\ds \frac{1}{2}\bSd:\dot{\bCd}-\rhoId\left(\dot{\psid}+\dot{T}\,\sd\right) - \frac{1}{T}\bQd\cdot\Gradd(T) \geq 0~.
\label{e:con_second_law_3D_non_polar}
\eqe
$\psid$ can be considered as a function of $\bCd_{\!\text{e}}$ and $T$ or $\bCd$ and $T$. First, the former one is considered,
\eqb{lll}
\psid \is \psid\left(\bCd_{\!\text{e}},T\right)~,
\eqe
so that its time derivative is
\eqb{lll}
\dot{\psid} \is \ds \pa{\psid}{\bCd_{\!\text{e}}}:\dot{\bCd}_{\!\text{e}}+\pa{\psid}{T}\dot{T}~.
\eqe
Using the time derivative of Eq.~(\ref{e:3D_right_cauchy_diff}.1), it can be shown that
\eqb{lll}
\dot{\bCd}_{\!\text{e}} \is \bFd^{-\text{T}}_{\!\text{T}}\,\dot{\bCd}\,\bFd^{-1}_{\!\text{T}} +\left(\bHd+\bHd^{\text{T}}\right)\,\dot{T}~,
\eqe
with
\eqb{lll}
\bHd \dis \bFd^{-\text{T}}_{\!\text{T}\,,T}\,\bCd\,\bFd^{-1}_{\!\text{T}}~.
\eqe
Here, $(\bullet)_{,T}$ indicates the partial derivative with respect to the absolute temperature.
Substituting these relations into Eq.~\eqref{e:con_second_law_3D_non_polar} results in
\eqb{lll}
\ds \left(\frac{1}{2}\bSd- \bFd^{-1}_{\!\text{T}}\rhoId \pa{\psid}{\bCd_{\!\text{e}}}\bFd^{-\text{T}}_{\!\text{T}}\right) :\dot{\bCd}- \rhoId\left(\pa{\psid}{\bCd_{\!\text{e}}}:\left(\bHd+\bHd^{\text{T}}\right) +\pa{\psid}{T}+\sd\right)\dot{T} - \frac{1}{T}\bQd\cdot\Gradd(T) \geq 0~.
\eqe
The latter relation should be positive for all admissible states of the continua. So the final relations for $\bSd$ and $\sd$ can be obtained as
\eqb{lll}
\ds
\ds \bSd \is \ds 2\rhoId\,\bFd_{\text{T}}^{-1}\,\pa{\psid}{\bCd_{\!\text{e}}}\,\bFd_{\text{T}}^{-\text{T}}~.
\label{e:3D_nonpolar_S}
\eqe
and
\eqb{lll}
\ds \sd\is \ds -\pa{\psid}{T}-\pa{\psid}{\bCd_{\!\text{e}}}:\left(\bHd+ \bHd^{\text{T}}\right)~.
\label{e:3D_nonpolar_s}
\eqe
Furthermore, a Mandel type stress in the intermediate configuration can be obtained by pushing forward (see Eq.~\eqref{e:push_forward}) $\bSd$ by $\bFd_{\text{T}}$ (denoted by $\rhd(\bF_{\text{T}})$) as
\eqb{lll}
\ds
\ds \bSd_{\text{T}} \dis \ds \frac{1}{\JTd}\,\bS^{\sharp\rhd(\bF_{\!\text{T}})} = \frac{1}{\JTd}\,\bFd_{\!\text{T}}\,\bSd\,\bFd^{\text{T}}_{\!\text{T}} =\frac{2\rhoId}{\JTd}\,\pa{\psid}{\bCd_{\!\text{e}}}~,
\eqe
where $\JTd$ is
\eqb{lll}
\JTd \is \detd\,\bFd_{\!\text{T}}~.
\eqe
Using the obtained relations for $\bSd$ and $\sd$, Eq.~\eqref{e:con_laws_3D_non_polar} can then be simplified as
\eqb{lll}
\ds \rhoId\,T\,\dot{\sd}=\rhoId\,\rd-\Divd\,\bQd~.
\eqe
If $\psid=\psid\left(\bCd,T\right)$, another set of relations can be obtained for $\bSd$ and $\sd$, i.e.~\citep{Holzapfel2000_01}
\eqb{lll}
\ds \bSd \is \ds \left(2\rhoId\,\pa{\psid}{\bCd}\right)_{\dot{T}=0},\\[3mm]
s \is \ds -\left(\pa{\psid}{T}\right)_{\dot{\bCd}=0}~.
\eqe
These relations are more simpler than the previous ones. However, the zero stress condition under pure thermal deformation should be satisfied by $\psid$. For $\psid=\psid\left(\bCd,T\right)$, the zero stress condition causes difficulties in the development of material models, and isothermal functional forms should be modified to resolve this issue. But it is straight forward for the case of $\psid(\bCd_{\!\text{e}},T)$, where the classical functional can be used (by replacement of $\bCd$ by $\bCd_{\!\text{e}}$). The temperature dependence is for example included by considering material parameters to be temperature dependent.
\subsection{Weak forms of 3D non-polar continua}
The weak form of the equilibrium relation can be obtained by its multiplication of an admissible variation $\delta \bx$ and integration over the domain as
\eqb{lll}
\ds  \int_{\mathcal{B}}{\rhod\,\delta\bxd\cdot \dot{\bvd}\,\dif v}\is \ds \int_{\mathcal{B}}{ \delta\bxd\cdot \divd\,\bsigd\,\dif v}+\int_{\mathcal{B}}{ \rhod\,\delta\bxd\cdot \bff\,\dif v}~,~\forall~\delta \bxd \in \mathcal{V}~,
\label{e:con_mon_weak_cur_01}
\eqe
where $\mathcal{V}$ is the space of admissible variations.
Eq.~\eqref{e:con_mon_weak_cur_01} can be simplified by the divergence theorem as
\eqb{lll}
\ds \int_{\mathcal{B}}{\gradd\,(\delta\bvd):\bsigd\,\dif v}+\int_{\mathcal{B}}{ \rhod\,\delta\bxd\cdot \dot{\bvd}\,\dif v}\is \ds \int_{\partial_{\bt}\mathcal{B}}{ \delta\bxd\cdot\btd\,\dif s}+\int_{\mathcal{B}}{ \rhod\,\delta\bxd\cdot \bff\,\dif v} ~,~\forall~\delta \bxd \in \mathcal{V}~.
\label{e:con_mon_weak_cur_02}
\eqe
Similarly, by multiplication of an admissible variation $\delta\theta$ and use of the divergence theorem, the weak form of the energy equilibrium can be written as
\eqb{lll}
\ds \int_{\mathcal{B}}{ \delta\theta\,\rhod\,T\,\dot{\sd}\,\dif v} \is \ds \int_{\mathcal{B}}{ \gradd(\delta\theta)\cdot\bqd\,\dif v}+\int_{\mathcal{B}}{ \delta\theta\,\rhod\,r\,\dif v}-\int_{\partial_{\bq}\mathcal{B}}{ \delta\theta\,\bqd\cdot\bnu\,\dif s}~,~\forall~\delta \theta \in \mathcal{V}~.
\eqe
\section{Surface continua}\label{s:Surf_shell}
In this section, the kinematics and balance laws for surface continua are obtained as a special case of 3D continua.
\subsection{Curvilinear description of deforming surfaces}
\begin{figure}[h]
    \begin{subfigure}[t]{1\textwidth}
        \centering
    \includegraphics[height=100mm]{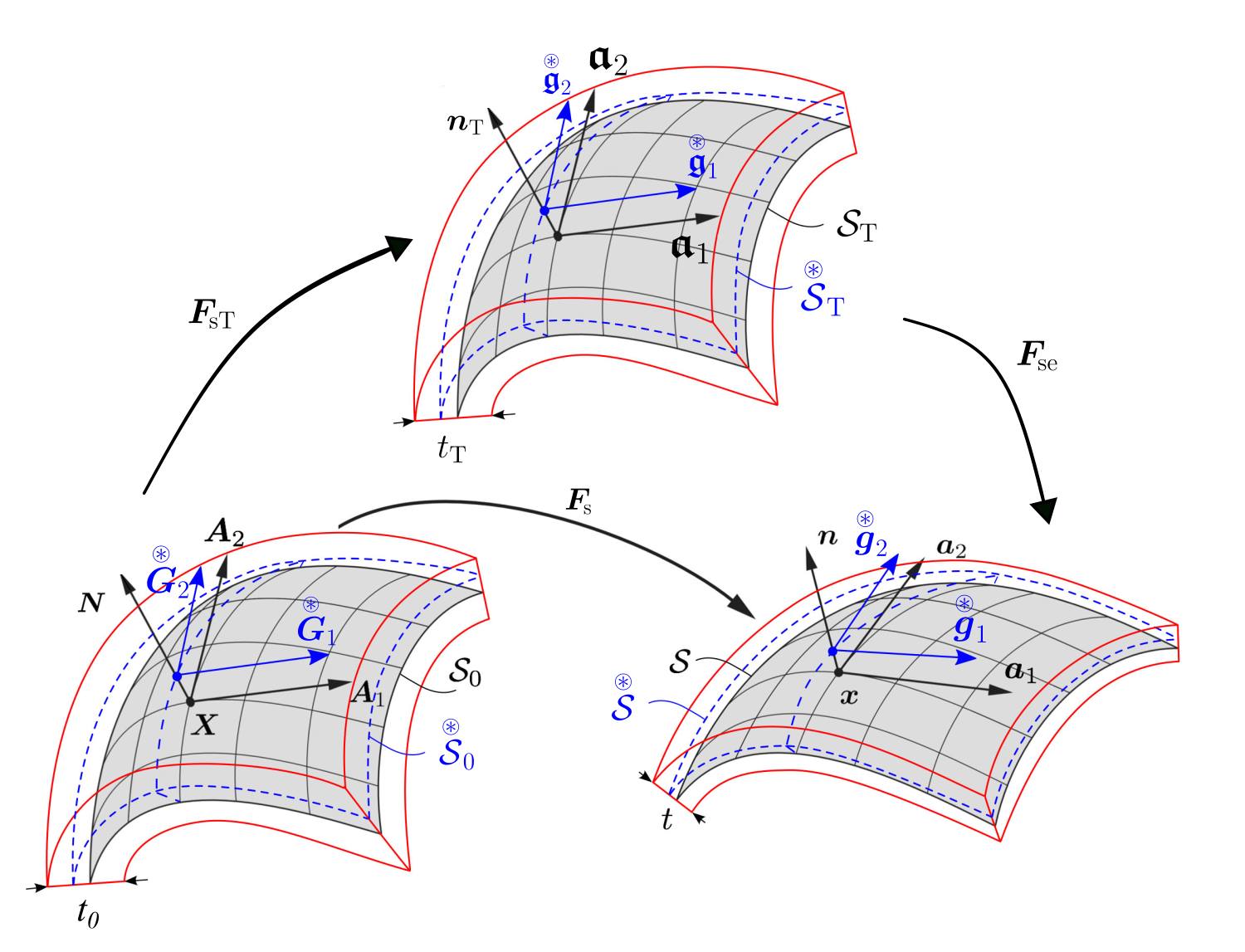}
    \end{subfigure}
    \caption{Multiplicative deformation decomposition of layered surface continua. $\mcalS_{0}$, $\mcalSTz$ and $\mcalS$ denote the mid-surface in reference, intermediate and current configuration and their corresponding tangent vectors are $\bA_{i}$, $\baTzui{i}$ and $\ba_{i}$. $\mcalScazF_{0}$, $\mcalScaTzF$ and $\mcalScazF$  are the reference, intermediate and current configuration of shell layers (through the thickness) and their corresponding tangent vectors are $\bgcazF_{i}$, $\bgcaTzF_{i}$ and $\bGcazF_{i}$. $\bFz$, $\bFTz$ and $\bFez$ are total, thermal and elastic surface deformation gradients.\label{f:Surface_thermal_kin}}
\end{figure}

The parametric description of the mid-surface in the reference and the current configuration can be written as (see Fig.~\ref{f:Surface_thermal_kin})
\eqb{lllll}
\bXz \is \ds \bXz(\xi^{\alpha})~,\quad \bxz \is \ds \bxz(\xi^{\alpha},t)~,
\eqe
where $\xi^{\alpha}$ are convective coordinates and the corresponding tangent vectors are
\eqb{lllll}
\bA_{\alpha} \is \ds \pa{\bXz}{\xi^{\alpha}}~,\quad \ba_{\alpha} \is \ds \pa{\bxz}{\xi^{\alpha}}~.
\eqe
Next, the co-variant metric is obtained by the inner product, i.e.
\eqb{lllll}
\Auab \is \bA_{\alpha}\cdot\bA_{\beta}~,\quad \auab \is \ba_{\alpha}\cdot\ba_{\beta}~.
\eqe
The contra-variant metric is the inverse of the co-variant metric, i.e.
\eqb{lllll}
[\Aab] \is [\Auab]^{-1}~,\quad [\aab] \is [\auab]^{-1}~.
\eqe
The dual tangent vectors can then be defined as
\eqb{lllll}
\bA^{\alpha} \dis \Aab\,\bA_{\beta}~,\quad \ba^{\alpha} \dis \aab\,\ba_{\beta}~.
\eqe
The unit normal vector of the surface in its reference and current configuration can then be written as
\eqb{lllll}
\bN \is \ds \frac{\bA_{1}\times\bA_{2}}{\|\bA_{1}\times\bA_{2}\|}~,\quad \bn \is \ds \frac{\ba_{1}\times\ba_{2}}{\|\ba_{1}\times\ba_{2}\|}~.
\eqe
The 3D identity tensor $\bolds{1}$ can be written as
\eqb{l}
\bolds{1} = \bI+\bN\otimes\bN = \bi+\bn\otimes\bn~,
\eqe
where $\bI$ and $\bi$ are the surface identity tensors in the reference and current configuration that can be written as
\eqb{lll}
\bI \is \Auab\,\bA^{\alpha}\otimes\bA^{\beta} = \Aab\,\bA_{\alpha}\otimes\bA_{\beta}=\bA^{\alpha}\otimes\bA_{\alpha}=\bA_{\alpha}\otimes\bA^{\alpha}~,
\eqe
\eqb{lll}
\bi \is  \auab\,\ba^{\alpha}\otimes\ba^{\beta} = \aab\,\ba_{\alpha}\otimes\ba_{\beta}=\ba^{\alpha}\otimes\ba_{\alpha}=\ba_{\alpha}\otimes\ba^{\alpha}~.
\eqe
The surface gradient operator for a general scalar, surface vector and surface tensor\footnote{Surface vectors and tensors have only in-plane components.} in the reference and current configuration are
\eqb{lll}
\ds \Gradz\, \Phi \dis \ds \Phi_{,\alpha}\,\bA^{\alpha}~,\\[3mm]
\ds \Gradz\, \buz \dis \ds \buzi{\,,\beta}\otimes\bA^{\beta} = u^{\alpha}_{;;\beta}\,\bA_{\alpha}\otimes\bA^{\beta}~,\\[3mm]
\ds \Gradz\, \bUz \dis \ds \bUzi{\,,\gamma}\otimes\bA^{\gamma} =U^{\alpha\beta}_{~~;;\gamma}\,\bA_{\alpha}\otimes\bA_{\beta}\otimes\bA^{\gamma}~\\
\eqe
and
\eqb{lll}
\ds \gradz\, \Phi \dis \ds \Phi_{,\alpha}\,\ba^{\alpha}~,\\[3mm]
\ds \gradz\, \bvz \dis \ds \bvzi{\,,\beta}\otimes\ba^{\beta}= v^{\alpha}_{~;\beta}\,\ba_{\alpha}\otimes\ba^{\beta}~,\\[3mm]
\ds \gradz\, \bVz \dis \ds \bVzi{\,,\gamma}\otimes\ba^{\gamma} = V^{\alpha\beta}_{~~~;\gamma}\,\ba_{\alpha}\otimes\ba_{\beta}\otimes\ba^{\gamma}~,\\
\eqe
where ``;;'' and ``;'' are the co-variant derivatives in the reference and current configuration. They are defined as \citep{Itskov2015_01}
\eqb{lll}
u^{\alpha}_{;;\beta} \dis u^{\alpha}_{,\beta}+u^{\gamma}\,\GammaIz{\gamma}{\beta}{\alpha}~,\\[3mm]
\ds U^{\alpha\beta}_{~~;;\gamma} \dis \ds U^{\alpha\beta}_{~~,\gamma}+U^{\delta\beta}\,\GammaIz{\delta}{\gamma}{\alpha}+U^{\alpha\delta}\,\GammaIz{\delta}{\gamma}{\beta}~,\\[3mm]
v^{\alpha}_{;\beta} \dis v^{\alpha}_{,\beta}+v^{\gamma}\,\Gammaz{\gamma}{\beta}{\alpha}~,\\[3mm]
V^{\alpha\beta}_{~~~;\gamma} \dis V^{\alpha\beta}_{~~~,\gamma}+V^{\delta\beta}\,\Gammaz{\delta}{\gamma}{\alpha}+V^{\alpha\delta}\,\Gammaz{\delta}{\gamma}{\beta}~,
\eqe
where the surface Christoffel symbols, $\GammaIz{\alpha}{\beta}{\gamma}$ and $\Gammaz{\alpha}{\beta}{\gamma}$, are defined as
\eqb{lllll}
\GammaIz{\alpha}{\beta}{\gamma} \dis \bA_{\alpha,\beta}\cdot\bA^{\gamma}~, \quad \Gammaz{\alpha}{\beta}{\gamma} \dis \ba_{\alpha,\beta}\cdot\ba^{\gamma}~.
\eqe
The co-variant derivative of the surface tangent and dual tangent vector are defined as \citep{Itskov2015_01}
\eqb{llllll}
\bA_{\alpha;;\beta} \dis \bA_{\alpha,\beta}-\bA_{\gamma}\,\GammaIz{\alpha}{\beta}{\gamma}~,\quad & \bA^{\alpha}_{;;\beta} \dis \bA^{\alpha}_{,\beta}+\bA^{\gamma}\,\GammaIz{\gamma}{\beta}{\alpha}~,\\[3mm]
\ba_{\alpha;\beta} \dis \ba_{\alpha,\beta}-\ba_{\gamma}\,\Gammaz{\alpha}{\beta}{\gamma}~,\quad & \ba^{\alpha}_{;\beta} \dis \ba^{\alpha}_{,\beta}+\ba^{\gamma}\,\Gammaz{\gamma}{\beta}{\alpha}~.
\eqe
Similar to the 3D case, see Eqs.~\eqref{e:3D_Ricci_ref} and \eqref{e:3D_Ricci_cur}, the co-variant derivative of the metric tensor in the reference and current configuration are zero, i.e.~\citep{Sauer2017_01}
\eqb{l}
A_{\alpha\beta;;\gamma} = 0~,\quad  A^{\alpha\beta}_{~~;;\gamma}= 0~,\quad a_{\alpha\beta;\gamma} = 0~,\quad a^{\alpha\beta}_{~~;\gamma}= 0~,
\eqe
so that the contra-variant derivatives follow as
\eqb{lllll}
U_{\alpha\beta;;}^{~~~\gamma} \is A_{\alpha \delta}\,A_{\beta \xi}\,A^{\gamma \eta}\,U^{\delta\xi}_{~~;;\eta}~,\quad
V_{\alpha\beta;}^{~~\gamma} \is a_{\alpha \delta}\,a_{\beta \xi}\,a^{\gamma \eta}\,V^{\delta\xi}_{~~;\eta}~.
\eqe
Furthermore, the surface divergence of a surface vector and a surface second order tensor is defined as
\eqb{lllllllllll}
\Divz\, \buz \dis u^{\alpha}_{;;\alpha}~,\quad & \divz\, \bvz \dis v^{\alpha}_{;\alpha}~,\quad
\Divz\, \bUz \dis U^{\alpha\beta}_{~~;;\beta}\,\ba_{\alpha}~,\quad & \divz\, \bVz \dis V^{\alpha\beta}_{~~;\beta}\,\ba_{\alpha}~.\\[3mm]
\eqe
Using Ricci's theorem (Eqs.~\eqref{e:3D_Ricci_ref} and \eqref{e:3D_Ricci_cur}), it can be shown that
\eqb{lllll}
\bA_{\alpha\|\beta} \is \bA_{\alpha;;\beta}-\bN\,\Gamma^{3}_{0~\alpha\beta}=0~,\quad \ba_{\alpha|\beta} \is \ba_{\alpha;\beta}-\bn\,\Gamma^{3}_{\alpha\beta}=0~,
\eqe
where the co-variant components of the curvature tensor are defined by
\eqb{lllll}
\bIuab \dis \Gamma^{3}_{0~\alpha\beta}:=\bN\cdot\bA_{\alpha,\beta}=\bN\cdot\bA_{\alpha;;\beta}~,\quad \buab \dis \Gamma^{3}_{\alpha\beta}:=\bn\cdot\ba_{\alpha,\beta}=\bn\cdot\ba_{\alpha;\beta}~,
\eqe
and Gauss' formula can be written as
\eqb{lllll}
\bA_{\alpha;;\beta} \is \bIuab\,\bN~,\quad \ba_{\alpha;\beta} \is \buab\,\bn~.
\eqe
Using differentiation of $\bN\cdot\vAua=0$, $\bN\cdot\bN=1$ and their counterparts in the current configuration, it can be shown that
\eqb{lllllllll}
\Gamma^{3}_{0\,3\beta} \is 0~,\quad \Gamma^{3}_{3\beta} \is 0~,\quad \Gamma^{\beta}_{0\,3\alpha} \is -b^{\beta}_{0~\alpha}~,\quad \Gamma^{\beta}_{3\alpha} \is -b^{\beta}_{\alpha}~.
\eqe
Using these relations, Weingarten's formula can be obtained as
\eqb{lllll}
\bN_{\!,\alpha} \is -b^{\beta}_{0\,\alpha}\,\vAub~,\quad \bn_{,\alpha} \is -b^{\beta}_{\alpha}\,\vaub~,
\eqe
or based on the gradient operator as
\eqb{lllll}
\ds -\Gradz\,\bN = \bIuab\,\bA^{\alpha}\otimes\bA^{\beta} \dis \bb_{0}~,\quad \ds -\gradz\,\bn = \buab\,\ba^{\alpha}\otimes\ba^{\beta}\dis \bb~,
\eqe
which allows us to identify the curvature tensors as the surface gradient of the surface normals.
The contra-variant components of $\bb_{0}$ and $\bb$ can be obtained as
\eqb{lllll}
\bab_{0} \is \Aag\,b_{0\,\gamma\delta}\,\Adb~,\quad \bab \is \aag\,b_{\gamma\delta}\,\adb~.
\eqe
The mean and Gaussian curvatures $H_{0}$ and $\kappa_{0}$ ($H$ and $\kappa$) of the reference (current) configuration are defined as
\eqb{l}
H_{0} := \ds  \frac{1}{2}\bb_{0}:\bI = \frac{1}{2}b^{\alpha}_{0\,\alpha}~,~H := \ds  \frac{1}{2}\bb:\bi = \frac{1}{2}b^{\alpha}_{\alpha}~,
\label{e:surf_mean_gauss1}
\eqe
\eqb{l}
\kappa_{0} := \ds \det(\bb_{0}) =\frac{\det[\bIuab]}{\det[\Augd]}~,~\kappa := \ds \det(\bb) =\frac{\det[\buab]}{\det[\augd]}~,
\eqe
where the surface determinant for a general surface tensor $\bTz=T_{\alpha\beta}\,\bl^{\alpha}\otimes\bk^{\beta}$ is defined as\footnote{$\bl^{\alpha}$ and $\bk^{\beta}$ are general base vectors, e.g.~$\ba^{\alpha}$ and $\bA^{\alpha}$.}
\eqb{lll}
\ds \detz\,\bTz \dis \ds \frac{[(\bTz \cdot\by_{1})\times(\bTz \cdot\by_{2})]\cdot\by_{4}}{[\by_{1}\times\by_{2}]\cdot\by_{3}}~,
\eqe
where $\by_{1}$ and $\by_{2}$ lie in the surface spanned by $\bk^{\beta}$, i.e.~they do not have any out of plane components and are not parallel to each other, and $\by_{3}$ is the unit normal vector to the surface. $\by_{4}$ is the unit normal vector to the plane spanned by $\bl^{\alpha}$. A similar definition of the determinant of surface tensors can be found in \citet{Javili2014_01}. Alternative to Eq.~\eqref{e:surf_mean_gauss1}, the mean and Gaussian curvatures can be defined from the principal curvatures $\kappa_{\alpha}$ as
\eqb{l}
H_{0} = \ds \frac{1}{2}(\kappa_{0\,1}+\kappa_{0\,2})~,~H = \ds \frac{1}{2}(\kappa_1+\kappa_2)~,
\eqe
\eqb{l}
\kappa_{0} = \ds \kappa_{0\,1}\,\kappa_{0\,2}~,~\kappa = \ds \kappa_1\,\kappa_2~,
\eqe
where $\kappa_{0\,\alpha}$ and $\kappa_{\alpha}$ are the eigenvalues of matrix $[\bIuab\,\Abg]$ and $[\buab\,\abg]$, respectively.\\
Assuming Kirchhoff-Love kinematics, the position of a point on the shell layers $\mcalScaz_{0}$ and $\mcalScaz$ can be written as
\eqb{lllll}
\bXcaz \dis \bXz+\xi\,\bN~,\quad \bxcaz \dis \bxz+\xi\,\lambda_{3}\,\bn~,
\eqe
where $\xi$ is the thickness coordinate and $\lambda_{3}$ is the stretch in the direction of the surface normal. The tangent vectors of $\bXcaz$ and $\bxcaz$ follow as
\eqb{lll}
\ds \bGcaz_{\alpha} \dis \ds \pa{\bXcaz}{\xi^{\alpha}} =\vAua-\xi\,b_{0\,\alpha}^{\gamma}\,\vAug=\vAua-\xi\,b_{0\,\alpha\gamma}\,\vAg~,
\eqe
\eqb{lll}
\ds \bgcaz_{\alpha} \dis \ds \pa{\bxcaz}{\xi^{\alpha}} =\vaua-\xi\,b_{\alpha}^{\gamma}\,\vaug=\vaua-\xi\,\lambda_{3}\,b_{\alpha\gamma}\,\vag~,
\eqe
and their corresponding metric up to leading order in $\xi$ are\footnote{Assuming $\lambda_{3}\simeq 1$.}
\eqb{lllll}
\bGcaz_{\alpha\beta} \dis \bGcaz_{\alpha}\cdot\bGcaz_{\beta} \simeq \Auab -2\xi\,\bIuab~,\quad \bgcaz_{\alpha\beta} \dis \bgcaz_{\alpha}\cdot\bgcaz_{\beta} \simeq \auab -2\xi\,\buab~.
\label{e:surf_metic_star}
\eqe
\subsection{Kinematics of surface deformation}\label{s:Surf_shell_kin_deform}
The deformation gradient of layer $\mcalScaz$ follows from Eq.~\eqref{3D_kin_def_grad} as
\eqb{lll}
\bFcaz \is \ds \pa{\bxcaz}{\bXcaz}=\bFaz+\lambda_{3}\,\bn\otimes\bN~,
\eqe
and it transforms $\dif\bXcaz$ to $\dif\bxcaz$ as
\eqb{lll}
\dif\bxcaz \is \ds \bFcaz\,\dif\bXcaz~.
\eqe
The in-plane part of $\bFcaz$, $\bFaz$, can be written based on the tangent vectors as
\eqb{lll}
\bFaz \is \ds \pa{\bxcaz}{\xi^{\alpha}}\otimes\pa{\xi^{\alpha}}{\bXcaz} =\bgcaz_{\alpha}\otimes~\bGcaz^{\alpha}~,
\eqe
and it transforms $\bGcaz_{\alpha}$ as
\eqb{lll}
\bgcaz_{\alpha} \is \ds \bFcaz\,\bGcaz_{\alpha}= \bFaz\,\bGcaz_{\alpha}~.
\eqe
The mid-surface deformation gradient (i.e.~for layer $\sS$) is defined as
\eqb{lll}
\ds \bFIz \dis \ds \bFcaz(\xi=0)=\bFz +\lambda_{3}\,\bn\otimes\bN~,
\eqe
where $\bFz=\bFaz(\xi=0)$ is the in-plane deformation gradient of the mid-surface, i.e.
\eqb{lll}
\bFz \is \ba_{\alpha}\otimes\bA^{\alpha}~.
\eqe
The deformation gradient $\bFcaz$ can be decomposed into thermal and elastic parts $\bFcaTz$ and $\bFcaez$ as
\eqb{lll}
\ds \bFcaz \is \ds \bFcaez\,\bFcaTz~.
\eqe
$\bFcaTz$ pushes forward $\dif\bXcaz$ to the intermediate configuration as
\eqb{lll}
\dif\bXcaTz \is \bFcaTz\,\dif\bXcaz~.
\eqe
The tangent vectors in the intermediate configuration $\bgcaTzui{\alpha}$ are
\eqb{lll}
\bgcaTzui{\alpha} \is \bXcaTzi{,\alpha}~.
\eqe
$\bFcaez$ and $\bFcaTz$ are
\eqb{lll}
\bFcaez \is \bFaez +\lambda_{\text{e}\,3}\,\bn\otimes\bn_{\text{T}}~,
\eqe
\eqb{lll}
\bFcaTz \is \bFaTz +\lambda_{\text{T}\,3}\,\bn_{\text{T}}\otimes\bN~,
\eqe
where $\lambda_{\text{e}\,3}$ and $\lambda_{\text{T}\,3}$ are the elastic and thermal parts of the stretch in the thickness direction, and $\bn_{\text{T}}$ is the unit normal vector of the surface in the intermediate configuration.
The surface right Cauchy-Green deformation tensor $\bCcaz$ and its thermal and elastic part $\bCcaTz$ and $\bCcaez$ are
\eqb{lll}
\bCcaz \is \bFcaz^{\text{T}}\,\bFcaz =\bFcaTz^{\text{T}}\,\bCcaez\,\bFcaTz~,\\[3mm]
\bCcaTz \is \bFcaTz^{\text{T}}\,\bFcaTz= \bCaTz+\lambda_{\text{T}\,3}^{2}\,\bN\otimes\bN~,\\[3mm]
\bCcaez \is \bFcaez^{\text{T}}\,\bFcaez =\bCaez+\lambda_{\text{e}\,3}^2\,\bn_{\text{T}}\otimes\bn_{\text{T}}~,
\label{e:Surf_right_cauchy_diff}
\eqe
where $\bCcaTz$ and $\bCcaez$ are
\eqb{lll}
\bCaTz \is \bFaTz^{\text{T}}\,\bFaTz =\gcaTzi{\alpha\beta}\bGcaz^{\alpha}\otimes\bGcaz^{\beta}~,
\eqe
\eqb{lll}
\bCaez \is \bFaez^{\text{T}}\,\bFaez =\gcaz_{\alpha\beta}\,\bgcaTz^{\alpha}\otimes\bgcaTz^{\beta}~,
\eqe
with $\gcaTzi{\alpha\beta}=\bgcaTzui{\alpha}\cdot\bgcaTzui{\beta}$.
The velocity gradient $\blcaz$ can be written as
\eqb{lll}
\ds \blcaz \dis \ds \dot{\bFcaz}\,\bFcaz^{-1}= \dot{\gcaz}_{\alpha}\otimes\bgcaz^{\alpha} +\dot{\bn}\otimes\bn+\frac{\dot{\lambda}_{3}}{\lambda_{3}}\,\bn\otimes\bn~.
\eqe
The velocity gradient of the mid-surface is $\blIz=\blcaz(\xi=0)$. For the mid-surface in the intermediate configuration, the tangent vectors can be written as\footnote{$\bFTz$ and $\bFez$ are such that $\bFz=\bFez\,\bFTz$.}
\eqb{lll}
\baTzui{\alpha}\is \bFTz\,\bA_{\alpha}~,
\eqe
and the curvature tensor as
\eqb{lll}
\bbTz \is -\text{grad}_{\text{T}}\,\bn_{\text{T}}= \bTzui{\alpha\beta}\,\baTzi{\alpha}\otimes\baTzi{\beta}~,
\eqe
where $\bTzui{\alpha\beta}$ is
\eqb{lllll}
\bTzui{\alpha\beta} \is \baTzui{\alpha,\beta}\cdot\bn_{\text{T}}~\quad\text{with}\quad \baTzui{\alpha,\beta}\is \bFTzi{,\beta}\,\bA_{\alpha}+ \bFTz\,\bA_{\alpha,\beta}~.
\eqe
The curvature change relative to the intermediate configuration is
\eqb{lll}
\bkappa \is \bb^{\flat\lhd}-\bbTz^{\flat\lhd(\bF_{\text{T}})}=k_{\alpha\beta}\,\bA^{\alpha}\otimes\bA^{\beta}~,
\eqe
with
\eqb{lllllll}
\kappa_{\alpha\beta} \is b_{\alpha\beta}-\bTzui{\alpha\beta}~,\quad \bb^{\flat\lhd} \is \bFz^{\text{T}}\,\bb\,\bFz~,\quad \bbTz^{\flat\lhd(\bF_{\text{T}})} \is \bFTz^{\text{T}}\,\bbTz\,\bFTz~.
\eqe
The local area change between the reference and current configuration can is given by
\eqb{lll}
\Jz \dis \ds \pa{a}{A}=\detz\,\bFz~.
\eqe
\subsection{Surface stress and moment tensors}
The Cauchy stress tensor of the Kirchhoff-Love shell formulation, $\bsig_{\text{KL}}$, and the 3D Stress tensor $\bsigd$ can be connected by replacing $\bsigd$ by $\bsig_{\text{KL}}/t$, where $t$ is the current thickness of the shell\footnote{$\lambda_{3,i}$ is neglected in the current work (for example see Eq.~\eqref{e:surf_metic_star}). The contribution of this term is included in \citet{Cirak2001_01}. Furthermore, for Kirchhoff-Love shells, $\sigma^{33}=\sigma^{3\alpha}=0$. $\bsig_{\text{KL}}$ can be obtained from $\bsigd$ such that the energy of 3D and surface continua is equal, see \citet{Roohbakhshan2017_01}.}. $\bsig_{\text{KL}}$ (force per length) can be written as
\eqb{lll}
\bsig_{\text{KL}} \is \Nab\,\vaua\otimes\vaub+S^{\alpha}\,\vaua\otimes\bn~,
\label{e:KL_tensor}
\eqe
Similar to replacing $\bsigd$ by $\bsig_{\text{KL}}$, the moment vector $\bmd$ is replaced by $\bmz/t$. $\bmz$ can be written as
\eqb{lll}
\bmz \is \bmurz^{\text{T}}\,\bnu~,
\eqe
where $\bmurz$ has only in-plane components and can be written as
\eqb{lll}
\bmurz\is \murz^{\alpha\beta}\,\vaua\otimes\vaub~,
\eqe
so $\bmurz$ can be written based on any set of in-plane vectors that are non-parallel. It is common to write $\bmurz$ as \citet{Sauer2017_01}
\eqb{lll}
\bmurz^{\text{T}}=\bn\times\bmuz^{\text{T}}~,
\eqe
where $\bmuz$ is
\eqb{lll}
\bmuz \is -\Mab\,\vaua\otimes\vaub~,
\label{e:surf_mu_rotated_trans}
\eqe
which can be used to write $\bmurz$ as
\eqb{lll}
\bmurz^{\text{T}} \is \Mab\,(\vaub\times\bn)\otimes\vaua~.
\eqe
$\vaub\times\bn$ is perpendicular to $\bn$ and it can be expressed in terms of the tangent and normal vector to the boundary, $\btau_{\!\text{v}}$ and $\bnu$, as
\eqb{lll}
\vaub\times\bn \is \tau_{\beta}\,\bnu-\nu_{\beta}\,\btau_{\!\text{v}}~\quad \text{with}\quad \tau_{\beta} = \btau_{\!\text{v}}\cdot\vaub~;~\nu_{\beta} = \bnu\cdot\vaub~,
\label{e:surf_cross_aua_n}
\eqe
where the following relations has been used
\eqb{lll}
\vaub = \tau_{\beta}\,\btau_{\!\text{v}}+\nu_{\beta}\,\bnu~,
\eqe
\eqb{l}
\btau_{\!\text{v}}\times\bn = \bnu~;~\bnu\times\bn=-\btau_{\!\text{v}}~.
\eqe
Using Eq.~\eqref{e:surf_cross_aua_n}, it can be shown that
\eqb{lll}
\murz^{\alpha\beta}\,\vaub \is \Mab\,(\tau_{\beta}\,\bnu-\nu_{\beta}\,\btau_{\!\text{v}})~.
\label{e:surf_mu_bar_half_comp}
\eqe
$\murz^{\alpha\beta}$ can be obtained by multiplying both sides of Eq.~\eqref{e:surf_mu_bar_half_comp} by $\ba^{\gamma}$, giving
\eqb{lll}
\murz^{\alpha\beta} \is M^{\alpha\gamma}\,\left(\tau_{\,\gamma}\,\nu^{\beta}-\nu_{\gamma}\,\tau^{\beta}\right)~,
\eqe
where $\tau^{\beta} = \btau_{\!\text{v}}\cdot\vab$ and $\nu^{\beta} = \bnu\cdot\vab$.
Furthermore, $\bmz$ can be written in terms of $\bnu$ and $\btau_{\!\text{v}}$ as
\eqb{lll}
\bmz \is \bmurz^{\text{T}}\,\bnu=m_{\nu}\,\bnu+m_{\tau}\,\btau_{\!\text{v}}~,
\eqe
where $m_{\nu}$ and $m_{\tau}$ are
\eqb{lll}
m_{\nu} \is \Mab\,\nu_{\alpha}\,\tau_{\beta}~,\\
m_{\tau} \is -\Mab\,\nu_{\alpha}\,\nu_{\beta}~.
\eqe
\subsection{Conservation laws for surface continua}\label{s:Surf_shell_Cons_laws}
\citet{Sahu2017_01} and \citet{Sauer2018_01} develop a thermomechanical Kirchhoff-Love shell formulation directly from the surface balance laws. Here on the other hand, the conservation laws for surface continua are derived as a special case of their 3D counterparts.
\subsubsection{Mass balance}
Using Eq.~\eqref{e:3D_cons_mass_loc}, $\Jd=\lambda_{3}\,\Jz$, $t=t_{0}\,\lambda_{3}$, and
\eqb{lll}
\rhoIz \dis \rhoId\,t_{0}~,\\
\rhoz \dis \ds \rhod\,t=\frac{\rhoId}{\lambda_{3}\,\Jz}\,t_{0}\,\lambda_{3}= \frac{\rhoId}{\Jz}\,t_{0}=\frac{\rhoIz}{\Jz}~,
\eqe
the mass conservation law for a surface continua can be written as
\eqb{lll}
\rhoIz \is \rhoz\,\Jz~,
\label{e:surf_cons_mass_loc}
\eqe
where $\rhoIz$ and $\rhoz$ are the surface density in the reference and current configuration.
\subsubsection{Linear momentum balance}
The surface linear momentum balance can be obtained directly from 3D linear momentum balance. Substituting relation $\rhoz := \rhod\,t$ and replacing $\bsigd$ by $\bsig_{\text{KL}}/t$ within the 3D linear momentum balance, i.e.~Eq.~\eqref{e:con_mon_loc_cur}, the momentum balance for surface continua is\footnote{\label{fn:3D_2D_bodyforce}The body force $\bff$ is force per unit mass, so this quantity is the same for surface and volume continua, i.e.~$\bffz=\bffd$.}
\eqb{lll}
\divd\,\bsig^{\text{T}}_{\text{KL}}+\rhoz\,\bff \is \rhoz\,\dot{\bv}~.
\label{e:surf_mon_loc_form1}
\eqe
Using $\divd\,\bsig^{\text{T}}_{\text{KL}}=\bsig^{\text{T}}_{\text{KL},k}\cdot\ba^{k}
=\sigma^{ji}_{\text{KL}|j}\,\ba_{i}$ and $\ba_{3}=\ba^{3}=\bn$, $\divd\,\bsig^{\text{T}}_{\text{KL}}$ can be written
\eqb{lll}
\sigma^{ji}_{\text{KL}|j}\ba_{i}\is  \left[\sigma^{ji}_{\text{KL},j}+\sigma^{li}_{\text{KL}}\,\Gamma^{j}_{lj}+\sigma^{jl}_{\text{KL}}\,\Gamma^{i}_{lj}\right]\ba_{i}~.
\label{e:surf_mon_loc_form1_part1}
\eqe
Using $\sigma^{3\alpha}_{\text{KL}}=0$, $\ba_{3}\cdot\ba^{\gamma}=0$, $\Gamma^{3}_{3\alpha}=0$, $\Gamma^{3}_{\alpha\beta}=\buab$ and $\Gamma^{\alpha}_{3\beta}=-b^{\alpha}_{\beta}$, Eq.~\eqref{e:surf_mon_loc_form1_part1} can be simplified as
\eqb{lll}
\divd\,\bsig^{\text{T}}_{\text{KL}}\is \left[\sigma^{\beta\alpha}_{\text{KL}\,;\beta}+\sigma^{\beta 3}_{\text{KL}}\,\Gamma^{\alpha}_{3\beta}\right]\ba_{\alpha}
+ \left[\sigma^{\beta 3}_{\text{KL}\,;\beta}+\sigma^{\alpha\beta}_{\text{KL}}\,\Gamma^{3}_{\beta\alpha}\right]\ba_{3}\\[3mm]
\is \left[N^{\beta\alpha}_{~~;\beta}-b^{\alpha}_{\gamma}\,S^{\gamma}\right]\ba_{\alpha}
+ \bigm[S^{\gamma}_{~;\gamma}+b_{\alpha\gamma}\,N^{\alpha\gamma}\bigm]\bn~.
\eqe
\citet{Sauer2017_01} obtain this relation directly from the linear momentum balance for surfaces, while it is obtained here from the corresponding 3D balance law.
\subsubsection{Angular momentum balance}\label{s:Surf_shell_Cons_laws_angular}
Assuming $\bc=\bzero$, the angular momentum balance for surface continua requires\footnote{The body force couple $\bcd$ is force per unit mass, so this quantity is the same for surface and volume continua, i.e.~$\bcz=\bcd$.}
\eqb{lll}
\divd\,\bmurz^{\text{T}}+\bsEd:\bsig_{\text{KL}}\is \bzero~.
\eqe
Using Eqs.~\eqref{e:div_vec_cross_ten} and \eqref{e:surf_mu_rotated_trans}, it can be simplified into
\eqb{lll}
\vaua\times\left[\left(\Nab-b_{\gamma}^{\beta}\,M^{\gamma\alpha}\right)\vaub+\left(S^{\alpha}+M^{\beta\alpha}_{~~;\beta}\right)\bn\right]=\bzero~.
\eqe
This will always hold if and only if
\eqb{lll}
\sigab \is \Nab-b_{\gamma}^{\beta}\,M^{\gamma\alpha}
\label{e:surf_sigab}
\eqe
is symmetric and
\eqb{lll}
S^{\alpha}\is -M^{\beta\alpha}_{~~;\beta}~.
\label{e:surf_Sa}
\eqe
See \citet{Sauer2017_01} for an alternative proof.
\subsubsection{Energy balance}\label{s:Surf_shell_Cons_laws_energy}
Using the relations for $\blIz$ and $\bsEd$ given in Eqs.~\eqref{e:surf_velocity_grad} and \eqref{e:ang_spin_to_vel}, it can be shown that
\eqb{lllll}
\bwAz \dis \ds \frac{1}{2}\bsEd:\blIz^{\text{T}} &=& \ds \frac{1}{2}\left(w_{\alpha\beta}\,\vaa\times\vab-\dot{\bn}\times\bn+\bn\times\dot{\bn} +\frac{\dot{\lambda}_{3}}{\lambda_{3}}\,\bn\times\bn\right)\\[3mm]
&& &=& \ds \frac{1}{2}\left(w_{\alpha\beta}\,\vaa\times\vab+2\bn\times\dot{\bn}\right)~,
\eqe
where $w_{\alpha\beta}$ is given in Eq.~\eqref{e:surf_w_ab}. This relation and the following identity for the quadruple product
\eqb{lll}
(\bad\times\bbd)\cdot(\bcd\times\bdd)=(\bad\cdot\bcd)(\bbd\cdot\bdd)-(\bad\cdot\bdd)(\bbd\cdot\bcd)
\eqe
can be used to obtained the relation
\eqb{lll}
\bwAz\cdot\bmurz^{\text{T}} \is \bwAz\cdot\bn\times\bmuz^{\text{T}}\\[3mm]
\is \ds -(\dot{\bn}\cdot\vaub)\,M^{\alpha\beta}\vaua = \dot{\bn}\cdot\bmuz^{\text{T}}~.
\label{e:surf_ang_vel_moment}
\eqe
Using Eqs.~\eqref{e:pull_backward}, \eqref{e:KL_tensor}, \eqref{e:surf_sigab}, \eqref{e:surf_Sa}, and the relations for $w_{\alpha}$, $\blIz$, $\dot{a}_{\alpha\beta}$, $\dot{b}_{\alpha\beta}$ in Eqs.~\eqref{e:surf_w_a}, \eqref{e:surf_velocity_grad}, \eqref{e:auab_rate} and \eqref{e:buab_rate}, it can be shown that
\eqb{lll}
\bsig^{\text{T}}_{\text{KL}}:\blIz \is \Nab\,w_{\alpha\beta}-S^{\alpha}\dot{\bn}\cdot\vaua\\[3mm]
\is \sigab\,w_{\alpha\beta}+b^{\beta}_{\gamma}\,M^{\gamma\alpha}\,w_{\alpha\beta}-M^{\beta\alpha}_{~~;\beta}\,w_{\alpha}\\[3mm]
\is \ds \frac{1}{2}\sigab\,\left(w_{\alpha\beta}+ w_{\beta\alpha}\right)+M^{\mu\alpha}\left(w_{\alpha\beta}\,b^{\beta}_{\mu}+ w_{\alpha;\mu}\right)-\left(M^{\beta\alpha}\,w_{\alpha}\right)_{;\beta}\\[3mm]
\is \ds \frac{1}{2}\sigab\,\dot{a}_{\alpha\beta}+ M^{\beta\alpha}\,\dot{b}_{\alpha\beta}-\left(M^{\beta\alpha}\vaua\cdot\,w_{\lambda}\,\ba^{\lambda}\right)_{;\beta}\\[3mm]
\is \ds \frac{1}{2}\bsigz^{\sharp\lhd}:\dot{\bC}-\bmuz^{\sharp\lhd}:\dot{\bb}^{\flat\lhd}-\divz\left(\dot{\bn}\cdot\bmuz^{\text{T}}\right)~,
\label{e:KL_2_vel}
\eqe
where $\bsigz^{\sharp\lhd}$ and $\bmuz^{\sharp\lhd}$ are obtained by using the pull back operator (see Eq.~\eqref{e:pull_backward}) as
\eqb{lll}
\ds \bsigz^{\sharp\lhd} \is \ds \bFz^{-1}\,\bsigz\,\bFz^{-\text{T}} =\sigab\,\vAua\otimes\vAub~~\text{with}~\bsigz=\sigab\,\vaua\otimes\vaub~,
\eqe
\eqb{lll}
\ds \bmuz^{\sharp\lhd} \is \ds \bF^{-1}\,\bmuz\,\bF^{-\text{T}} =-M^{\alpha\beta}\,\vAua\otimes\vAub~~\text{with}~\bmuz =-\Mab\,\vaua\otimes\vaub~.
\eqe
The following relation for the surface energy balance can be obtained by substituting Eqs.~\eqref{e:KL_2_vel} and \eqref{e:surf_ang_vel_moment} into Eq.~\eqref{e:con_eng_loc_cur_form2}, as\footnote{The heat source $\rd$ is force per unit mass, so this quantity is the same for surface and volume continua, i.e.~$\rz=\rd$}
\eqb{lll}
\ds \rhoz\,\uzt-\rhoz\,r-\frac{1}{2}\,\bsigz^{\sharp\lhd}:\bCzt+ \bmuz^{\sharp\lhd}:\dot{\bb}^{\flat\lhd} +t\,\divd\,\accentset{\circledast}{\bq} \is 0~.
\label{e:surf_eng_loc_cur_form1}
\eqe
Using the relations for stress and bending tensor from \ref{s:Surf_shell_Cons_laws_second_law}, the energy balance can be written as
\eqb{lll}
\ds \rhoz\,T\,\szt=\rhoz\,r-t\,\divd\,\bqcaz~.
\label{e:surf_eng_loc_cur_form2}
\eqe
\subsubsection{The second law of thermodynamics}\label{s:Surf_shell_Cons_laws_second_law}
Using Eqs.~\eqref{e:KL_2_vel}, \eqref{e:surf_ang_vel_moment} and \eqref{e:con_second_law_loc_cur_form}, the second law of thermodynamics for surface continua can be written as
\eqb{lll}
\ds \rhoz\,T\,\szt-\rhoz\,\uzt+\frac{1}{2}\,\bsigz^{\sharp\lhd}:\bCzt- \bmuz^{\sharp\lhd}:\dot{\bb}^{\flat\lhd} -\frac{t}{T}\accentset{\circledast}{\bq}\cdot\gradd(T) \geq 0~.
\label{e:con_second_law_surf_cur}
\eqe
The surface Helmholtz energy per unit mass is defined as
\eqb{lll}
\psiz \dis \uz-T\,\sz~,
\eqe
and its time derivative is
\eqb{lll}
\dot{\psiz} \is \uzt-\dot{T}\,\sz-T\,\szt~.
\eqe
Using the previous relation, Eq.~\eqref{e:con_second_law_surf_cur} can be simplified as
\eqb{lll}
\ds \frac{1}{2}\,\bsigz^{\sharp\lhd}:\bCzt-\bmuz^{\sharp\lhd}:\dot{\bb}^{\flat\lhd} -\rhoz\left(\dot{\psiz}+\dot{T}\,\sz\right)-\frac{t}{T}\bqcaz\cdot\gradd(T) \geq 0~.
\eqe
It is assumed that $\psiz$ can be written as
\eqb{lll}
\psiz \is \psiz\left(\bCez, \bkappa,T\right)~,
\eqe
so that its time derivative is
\eqb{lll}
\psizt \is \ds \pa{\psi}{\bCez}:\bCezt+\pa{\psi}{\bkappa}:\dot{\bkappa}+\pa{\psi}{T}\dot{T}~.
\eqe
Using the time derivative of Eq.~(\ref{e:Surf_right_cauchy_diff}.1), it can be shown that
\eqb{lll}
\bCezt \is \bFTz^{-\text{T}}\,\bCzt\,\bFTz^{-1} +\left(\bHz+\bHz^{\text{T}}\right)\,\dot{T}~,
\eqe
with
\eqb{lll}
\bHz \is \bFTTz^{-\text{T}}\,\bCz\,\bFTz^{-1}~,
\eqe
and
\eqb{lll}
\dot{\bkappa} \is \dot{\bb}^{\flat\lhd}-\dot{\bbTz}^{\flat\lhd}~.
\eqe
Substituting these relations into Eq.~\eqref{e:con_second_law_surf_cur} results in
\eqb{lll}
\ds \left(\frac{1}{2}\bsigz^{\sharp\lhd}- \bFTz^{-1}\rhoz \pa{\psi}{\bCez}\bFTz^{-\text{T}}\right):\bCzt -\left(\bmuz^{\sharp\lhd}+\pa{\psiz}{\bk}\right):\dot{\bb}^{\flat\lhd}
\\[3mm]- \ds \rhoz\left(\pa{\psiz}{\bCez}:\left(\bHz+\bHz^{\text{T}}\right)-\pa{\psiz}{\bk}:\bbTz^{\flat\lhd}_{,T} +\pa{\psiz}{T}+\sz\right)\dot{T} - \frac{t}{T}\accentset{\circledast}{\bq}\cdot\gradd(T) \geq 0~.
\eqe
The final relation for $\bsigz^{\sharp\lhd}$, $\bmuz^{\sharp\lhd}$ and $\sz$ can be obtained as
\eqb{lll}
\ds \bsigz^{\sharp\lhd} \is \ds 2\bFTz^{-1}\rhoz \pa{\psiz}{\bCez}\bFTz^{-\text{T}}~,
\label{e:surf_sig_back}
\eqe
\eqb{lll}
\ds \bmuz^{\sharp\lhd}\is \ds -\pa{\psiz}{\bkappa}
\label{e:surf_mu_back}
\eqe
and
\eqb{lll}
\ds \sz \is \ds -\pa{\psiz}{\bCez}:\left(\bHz+\bHz^{\text{T}}\right)+\pa{\psiz}{\bkappa}:\bbTz^{\flat\lhd}_{,T} -\pa{\psiz}{T}~.
\label{e:surf_s}
\eqe
With the equations of Sec.~\ref{s:Surf_shell_Cons_laws}, the thermomechanical Kirchhoff-Love shell formulation of \citet{Sahu2017_01} is extended by including the multiplicative decomposition of the deformation gradient and the heat transfer through the shell thickness. Further, the shell formulation of \citet{Sahu2017_01} is a componentwise formulation, while the current work considers a tensorial formulation.
\subsection{Weak forms}\label{s:Surf_shell_Cons_laws_weak_form}
Following \citet{Sauer2017_01}, the weak form of the equilibrium equation can be obtained as
\eqb{lll}
\ds \int_{\mathcal{S}}{\rhoz\,\delta \bxz\cdot\bvzt\,\dif a}\is \ds \int_{\mathcal{S}}{\delta \bxz\cdot\divd\,\bsig^{\text{T}}_{\text{KL}}\,\dif a}+
\int_{\mathcal{S}}{\rhoz\,\delta \bxz\cdot\bb\,\dif a}~,~\forall~\delta \bx \in \mathcal{V}~.
\label{e:surf_mon_weak_cur_01}
\eqe
Using the divergence theorem and $\btz=\bsig^{\text{T}}_{\text{KL}}\,\bnu$, this can be rewritten as
\eqb{lll}
\ds \int_{\mathcal{S}}{\gradz\,(\delta \bxz):\bsig^{\text{T}}_{\text{KL}}\,\dif a}+\int_{\mathcal{S}}{\rhoz\,\delta \bxz\cdot\bvzt\,\dif a}\is \ds \int_{\partial_{\bt}\mathcal{S}}{\delta \bxz\cdot\btz\,\dif l}+
\int_{\mathcal{S}}{\rhoz\,\delta \bxz\cdot\bb\,\dif a}~,~\forall~\delta \bx \in \mathcal{V}~.
\label{e:surf_mon_weak_cur_02}
\eqe
It can be shown that
\eqb{lll}
\ds M^{\beta\alpha}_{~~;\beta}\,\ba_{\alpha}\cdot\delta\bn \is \ds (M^{\beta\alpha}\,\ba_{\alpha}\cdot\delta\bn)_{;\beta}-M^{\beta\alpha}\,(\ba_{\alpha}\cdot\delta\bn)_{;\beta}\\[3mm]
\is \ds -\divz(\delta\bn\cdot \bmuz^{\text{T}})-M^{\beta\alpha}\,[\ba_{\alpha;\beta}\cdot\delta\bn+\ba_{\alpha}\cdot\delta\bn_{;\beta}]\\[3mm]
\is \ds -\divz(\delta\bn\cdot \bmuz^{\text{T}})-M^{\beta\alpha}\,\ba_{\alpha}\cdot\delta\bn_{,\beta}
\eqe
where $\ba_{\alpha;\beta}\cdot\delta\bn=0$ is used. Using the latter relation, \eqref{e:KL_tensor}, \eqref{e:surf_sigab}, \eqref{e:surf_Sa}), Weingarten $\ba_{,\alpha}=-b^{\beta}_{\alpha}\,\ba_{\beta}$ and $\delta \ba_{\alpha}\cdot\bn=-\ba_{\alpha}\cdot\delta\bn$, the relation $\gradz\,(\delta \bx):\bsig^{\text{T}}_{\text{KL}}$ can be written as
\eqb{lll}
\ds \gradz\,(\delta \bxz):\bsig^{\text{T}}_{\text{KL}} \is \delta \bxzi{\,,\eta}\otimes\ba^{\eta}:(\Nab\,\ba_{\beta}\otimes\ba_{\alpha}+S^{\alpha}\,\bn\otimes\ba_{\alpha})= \Nab\,\ba_{\beta}\cdot\delta\ba_{\alpha}+S^{\alpha}\,\bn\cdot\delta\ba_{\alpha}\\[3mm]
\is \sigab\, \delta\ba_{\alpha}\cdot\ba_{\beta}+b^{\beta}_{\gamma}\,M^{\gamma\alpha}\,\delta\ba_{\alpha}\cdot \ba_{\beta}-M^{\beta\alpha}_{~~;\beta}\,\bn\cdot\delta\ba_{\alpha}\\[3mm]
\is \ds \frac{1}{2}\sigab\, \delta\auab-M^{\beta\alpha}\,\delta\ba_{\alpha}\cdot \bn_{,\beta}+M^{\beta\alpha}_{~~;\beta}\,\delta\bn\cdot\ba_{\alpha}\\[3mm]
\is \ds \frac{1}{2}\,\bsigz^{\sharp\lhd}:\delta{\bCz}-M^{\beta\alpha}\,(\delta\ba_{\alpha}\cdot \bn_{,\beta}+\,\ba_{\alpha}\cdot\delta\bn_{;\beta})-\divz(\delta\bn\cdot \bmuz^{\text{T}})
\\[3mm]
\is \ds \frac{1}{2}\,\bsigz^{\sharp\lhd}:\delta{\bCz}-\bmuz^{\sharp\lhd}:\delta{\bb}^{\flat\lhd}-\divz(\delta\bn\cdot \bmuz^{\text{T}})~,
\eqe
where $\delta\buab=\delta\ba_{\alpha}\cdot \bn_{,\beta}+\,\ba_{\alpha}\cdot\delta\bn_{;\beta}$ is used. Using the latter relation, the mechanical weak form can be written as
\eqb{lll}
G_{\text{in}}+G_{\text{int}}=G_{\text{ext}}~,
\eqe
with
\eqb{lll}
G_{\text{in}} \is \ds \int_{\mathcal{S}}{\rho\,\delta \bxz\cdot\bvzt\,\dif a}~,
\eqe
\eqb{lll}
G_{\text{int}} \is \ds \int_{\mathcal{S}}{\frac{1}{2}\,\bsigz^{\sharp\lhd}:\delta\bCz\,\dif a}- \int_{\mathcal{S}}{\bmuz^{\sharp\lhd}:\delta{\bb}^{\flat\lhd}\,\dif a}
\eqe
and
\eqb{lll}
G_{\text{ext}} \is \ds \int_{\mathcal{S}}{\rhoz\,\delta \bxz\cdot\bff\,\dif a}+\int_{\partial_{\bt}\mathcal{S}}{\delta \bxz\cdot\btz\,\dif l}+\int_{\partial_{\bm}\mathcal{S}}{\delta\bn\cdot \bmuz^{\text{T}}\cdot\bnu\,\dif l}~.
\eqe
Multiplying Eq.~\eqref{e:surf_eng_loc_cur_form2} by an admissible test function, $\delta\theta$, the weak form of the energy balance can be written as
\eqb{lll}
\ds \int_{\mathcal{S}}{ \delta\theta\,\rhoz\,T\,\szt\,\dif a} \is \ds \int_{\mathcal{B}}{ \gradd(\delta\theta)\cdot\accentset{\circledast}{\bq}\,\dif v}+\int_{\mathcal{S}}{ \delta\theta\,\rhoz\,r\,\dif a}-\int_{\partial_{\bq}\mathcal{B}}{ \delta\theta\,\accentset{\circledast}{\bq}\cdot\bnu\,\dif s}~,~\forall~\delta \theta \in \mathcal{V}~.
\label{e:surf_weak_first_law_form1}
\eqe
In \ref{s:surf_objects_linearization}, the tensorial form for the linearization of the surface objects $\Jz$, $\bCz^{-1}$, $H$, $\kappa$ and  $\bb^{\flat\lhd}$ ($\bab$) are given. They can be used to linearize the weak forms without considering any specific coordinate system and thus generalize the curvilinear FE formulation of \citet{Duong2016_01}. If out-of-plane heat transfer is neglected, Eq.~\eqref{e:surf_weak_first_law_form1} reduces to
\eqb{lll}
\ds \int_{\mcalS}{ \delta\theta\,\rhoz\,T\,\szt\,\dif a} \is \ds \int_{\mcalS}{ \gradz(\delta\theta)\cdot\bqz\,\dif s}+\int_{\mathcal{S}}{ \delta\theta\,\rhoz\,r\,\dif a}-\int_{\partial_{\bq}\mcalS}{ \delta\theta\,\bqz\cdot\bnu\,\dif l}~,~\forall~\delta \theta \in \mathcal{V}~.
\eqe
However, if the temperature variation through the thickness has a major effect, the shell either needs to be discretized across the thickness or the temperature variation needs to be approximated by polynomials, see \citet{Surana1991_01}.
\section{Constitutive laws}\label{s:constitutive_laws}
In this section, lattice and material symmetry are discussed. Then, the evolution of symmetry groups is considered. Next, the constitutive laws for the heat flux and stress are discussed. Finally, two examples of the Helmholtz free energy are given for 3D and shell continua.
\subsection{Lattice and material symmetry}
A symmetry group of a crystal lattice is the set of all operations that leave the lattice indistinguishable from its initial configuration \citep{Tadmor2011_01}. The material symmetry group (the symmetry for physical properties such as the strain energy density) can have more symmetry operators than the crystal lattice \citep{Neumann1885_01,Newnham2005_01}. The Neumann principle mentions that ``The symmetry elements of any physical property of a crystal must include the symmetry elements of the point group of the crystal'' \citep[p.~20]{Nye1985_01} (see \citep{Neumann1885_01} in German). For example, the symmetry group of cubic crystals does not apply enough constraints on their constitutive laws in order to make them isotropic. However, optical properties of cubic crystals are isotropic and it means that the optical properties have higher symmetry than the lattice \citep[p.~20]{Nye1985_01}. These symmetry operators can be combined into structural tensors and these structural tensors can be used in the development of constitutive laws \citep{Zheng_1993_01,Zheng_1994_01,Zheng_1994_02,Menzel2003_01}. The lattice structure does not change for different temperatures and the crystal types do not depend on the temperature, if phase transformations are excluded. So, the thermal expansion of lattices does not change the lattice structure and structural tensors of constitutive laws for crystals \citep[p.~106-107]{Nye1985_01}.
\subsection{Evolution of the symmetry group}
If the preferred directions of a material, i.e.~the directions of anisotropy, transform with the total deformation gradient, these directions are the material direction\footnote{I.e.~the anisotropic directions in the reference and current configuration, $\byd_{0}$ and $\byd$, are connected by $\byd=\bFd\,\byd_{0}$.}. For example, the fiber directions in composite materials are material directions. In crystal plasticity, the plastic deformation gradient does not change the preferred directions\footnote{Unless twining occurs \citep{Kalidindi1998_01,Kalidindi2004_01,Clayton2011_01}.}. Only the elastic deformation gradient changes the lattice direction. Hence their preferred directions are not material directions. Assuming the preferred direction are the material direction, \citet{Reese2003_01} and \citet{Reese2008_01} push forward the preferred direction and normalize it. Structural tensors in the intermediate configuration are needed for the development of material models. Structural tensors in the reference and intermediate configuration can be defined as
\eqb{lll}
\bLd_{0I} \is \byd_{0}\otimes\byd_{0}~,
\eqe
\eqb{lll}
\bLd_{\text{T}} \is \byd_{\text{T}}\otimes\byd_{\text{T}}~,
\eqe
where $\byd_{0}$ and $\byd_{\text{T}}$ are the normalized preferred directions, e.g.~the fiber directions in fiber-reinforced composite materials, in the reference and intermediate configuration and are connected by
\eqb{lll}
\byd_{\text{T}} \is \ds \frac{\bFd_{\!\text{T}}\,\byd_{0}}{\|\bFd_{\!\text{T}}\,\byd_{0}\|}~.
\eqe
\citet{Cuomo2015_01} use a similar approach without normalization. \citet{Sansour2007_01} investigate the co-variant, contra-variant and mixed forms of the push forward for $\bL_{0}$. These forms are defined as
\eqb{lllll}
\ds \bLd^{\sharp\rhd(\bFd_{\!\text{T}})}_{\text{T}} \dis \ds \bFd_{\!\text{T}} \,\bLd^{\sharp}_{0} \,\bFd^{\text{T}}_{\!\text{T}} \is L^{ij}_{0}\,\bgTdui{i}\otimes\bgTdui{j}~,\\[3mm]
\ds \bLd^{\flat\rhd(\bFd_{\!\text{T}})}_{\text{T}} \dis \ds \bFd^{-\text{T}}_{\!\text{T}}\,\bLd^{\flat}_{0}\,\bFd^{-1}_{\!\text{T}} \is L_{0\,ij}\,\bgTdi{i}\otimes\bgTdi{j}~,\\[3mm]
\ds \bLd^{\backslash\rhd(\bFd_{\!\text{T}})}_{\text{T}} \dis \ds \bFd_{\text{T}}\,\bLd^{\backslash}_{0}\,\bFd^{-1}_{\text{T}} \is L^{i}_{0\,j}\,\bgTdui{i}\otimes\bgTdi{j}~,
\eqe
where $L^{ij}_{0}$, $L_{0\,ij}$ and $L^{i}_{0\,j}$ are the co-variant, contra-variant and mixed components of $\bLd_{0}$.
\citet{Schroder2002_01,Eidel2003_01,Sansour2006_01,Sansour2008_01,Dean2017_01} and \citet{Dean2018_01} use structural tensors of the reference configuration. The angle between the preferred directions for composite materials can change. \citet{Gong2016_01} consider the change of angles between preferred directions for woven composites (see also \citet{Peng2013_01,Alsayednoor2017_01} ).
\subsection{Thermal expansion}
Thermal expansion of anisotropic materials can be experimentally measured \citep{Dudescu2006_01}. For example, $\bFd_{\!\text{T}}$ can be written as \citep{Meissonnier2001_01}
\eqb{lll}
\bFd_{\text{T}} \is \ds e^{\balphad(\theta-\theta_{0})}~,
\eqe
where $\theta_{0}$ is a reference temperature, and $\balphad$ is a second order tensor that is related to thermal expansion\footnote{For monoclinic crystals, $\balphad$ can be written in matrix form as $\alpha_{ij}=\left[
                                                       \begin{array}{ccc}
                                                         \alpha_{11} & 0 & \alpha_{31} \\
                                                         0 & \alpha_{22} & 0 \\
                                                         \alpha_{31} & 0 & \alpha_{33} \\
                                                       \end{array}
                                                     \right]~,$ where $\alpha_{ij}$ are the material constants \citep{Nye1985_01,Wang2005_01}.} and can depend on temperature.
                                                     If $\balphad$ is temperature independent, $\dot{\bFd}_{\text{T}}$ can be written as
\eqb{lll}
\dot{\bFd}_{\text{T}} \is \dot{\theta}\,\balphad\,\bFd_{\text{T}}=\dot{\theta}\,\bFd_{\text{T}}\,\balphad~,
\eqe
where the latter relation is obtained by the coaxiality of $\balphad$ and $\bFTd$ (see \citet{Basar2000_01} for coaxiality). \citet{Wang2005_01} obtain anisotropic thermal expansion of monoclinic potassium lutetium tungstate single crystals (see \citet{Nye1985_01} for the structure of thermal expansion for other crystals).
\subsection{Heat flux}
The thermal conductivity tensor $\bkd$ should contain the symmetry group of the lattice. The symmetric and skew symmetric parts of the heat conductivity tensor are denoted in the following as $\bkd_{\text{sym}}$ and $\bkd_{\text{skew}}$. There is strong evidence that the thermal conductivity should be a symmetric tensor, but this can not be mathematically proven. This is discussed in the following.\\
First, $\bkd_{\text{skew}}$ cancels out in the Clausius-Planck inequality and does not contribute to entropy generation. Using Fourier's law and Eq.~(\ref{e:con_second_law_loc_cur_form2}.2), it can be proven that $\bkd_{\text{skew}}$ does not play role in entropy generation, i.e.
\eqb{lll}
\gamma_{\text{con}} \is \ds \frac{1}{T^2}\bkd:\gradd\,T\otimes\gradd\,T=\frac{1}{T^2}\bkd_{\text{sym}}:\gradd\,T\otimes\gradd\,T~,
\eqe
and $\bkd_{\text{skew}}$ also cancels out in the first law of thermodynamics for homogenous materials, i.e.~$\divd\,\bqd$ in Eq.~(\ref{e:con_eng_loc_cur_form2}) can be written as
\eqb{lll}
\ds \divd\,\bqd \is \ds \bkd:\paqq{T}{\bx}{\bx}+\divd\,(\bkd)\cdot\gradd\,T= \bkd_{\text{sym}}:\paqq{T}{\bx}{\bx}~,
\eqe
where the symmetry of $\partial^{2}\,T/\partial\,\bx\,\partial\,\bx$ and homogeneity ($\divd\,(\bkd)=0$) are used.
\citet{Casimir1945_01} (see also \citet{Groot1954_01,Groot1954_02}) shows that $\divd\,(\bkd-\bkd^{\text{T}})=0$ by using Onsager's principle \citep{Onsager1931_01,Onsager1931_02}\footnote{Onsager's principle is not true if magnetic and Coriolis forces exist. Furthermore, based on Onsager's principle, the heat flow should be the time derivative of a thermodynamical state variable and this statement is not true. Onsager's principle is not correctly applied in the original work of \citet{Onsager1931_01,Onsager1931_02} but in the corrected work by \citet{Casimir1945_01}.}. So, $\bkd_{\text{skew}}$ does not effect the first law of thermodynamics even for nonhomogeneous materials. If the heat conductivity of vacuum is considered to be zero, then the heat conductivity tensor will be symmetric \citep{Casimir1945_01}.\\
Second, a pure $\bkd_{\text{skew}}$ results in a spiral heat flux (see \citet[p.~205-207]{Nye1985_01} and \citet{Powers2004_01}), but a spiral heat flux has not been observed experimentally \citep{Soret1893_01,Voigt1910_01,Voigt1903_01}. It should be emphasized that zero heat conductivity in vacuum can not be proven and it is an assumption\footnote{See \citet{Mazur1997_01,Truesdell1969_01,Verhas2004_01,Cimmelli2014_01} for a discussion on the validity of Onsager's principle in \citet{Onsager1931_01,Onsager1931_02} and \citet{Casimir1945_01}.}. The thermal heat conductivity tensor can vary with strain \citep{Hu2013_01}.
\subsubsection{Heat conduction laws}
In this section, different heat conduction laws are discussed. Fourier's law can be written as \citep{Fourier1822_01,Lienhard2013_01}
\eqb{lll}
\ds \bqd \dis \ds -\bkd\,\gradd\,T~.
\eqe
Thermomechanical formulations based on Fourier's law are coupled problems of hyperbolic and parabolic types. The influence of the temperature variation reaches the entire continua instantly and disturbances propagates with infinite speed \citep{Hetnarski2008_01}. This issue can be resolved by considering more advanced heat transformation laws and considering a finite speed for propagation of thermal disturbances. \citet{Lord1967_01} extend Fourier's law by including a relaxation time. The Lord-Shulman model can be written as
\eqb{lll}
\ds \bqd \dis \ds -\bCd_{\text{rt}}\,\accentset{\circ}{\bqd} -\bkd\,\gradd\,T~,
\eqe
where $\bCd_{\text{rt}}$ is a second order tensor related to the relaxation time and the thermal propagation speed. $\accentset{\circ}{\bqd}$ is a suitable objective rate of the heat flux. This rate can be considered as the material time derivative \citep{Christov2005_01}, or the convected differentiation \citep{Christov2009_01,Khayat2011_01} as suggested by \citet{Oldroyd1950_01} for time dependent material constitutive laws. The convected differentiation is defined as
\eqb{lll}
\ds \frac{\mathfrak{d}\bqd}{\mathfrak{d}t}  \dis \ds \pa{\bqd}{t}+\left(\bv\cdot\gradd\right)\bqd-\gradd(\bv)\cdot\bqd+\divd(\bv)\,\bqd~,
\eqe
\citet{Green1972_01} propose another model by introducing the thermodynamic temperature \citep{Suhubi1975_01}
\eqb{lll}
\ds T^{\diamond} \dis \ds T+t_{1}\left[\dot{T}+b_{1}\,(T-T_{0})\,\dot{T}+b_{2}\,\dot{T}^{2}\right]~,
\eqe
where $t_{1}$, $b_{1}$ and $b_{2}$ are material constants. $T$ is replaced by $T^{\diamond}$ in the balance laws and Helmholtz free energy. This model is based on an extension of works by \citet{Muller1971_01} and \citet{Green1972_02}. The Lord-Shulman model does not change the balance laws, so the calibration of the material parameters and the numerical implementation are easier than \citet{Green1972_02}. See \citet{Chandrasekharaiah1986_01} and \citet{Hetnarski2008_01} for a comparison of the models of \citet{Green1972_02} and \citet{Lord1967_01}.
\subsubsection{Heat radiation}
The heat flux due to radiation $\bqd_{\text{r}}$ can be written as \citep{Lienhard2013_01}
\eqb{lll}
\ds \bqd_{\text{r}}\cdot\bnu \dis \ds -\epsilon\,\sigma\,T^{4}_{\text{s}}+q_{\text{rec}}~,
\eqe
where $T_{\text{s}}$ is the surface temperature in Kelvin, $\epsilon$ is emissivity and $\sigma$ is Stefan-Boltzmann constant. $q_{\text{rec}}$ is the amount of heat flux that is received by radiation. It can be written as
\eqb{lll}
\ds q_{\text{rec}}\dis \ds F_{\text{geo}}\,\epsilon\,\sigma\,T^{4}_{\text{ref}}~,
\eqe
where $T_{\text{ref}}$ is a reference temperature and $F_{\text{geo}}$ is a coefficient related to the surface geometry. $F_{\text{geo}}$ becomes unity for two parallel infinite flat plates or an enclosed continuum by another continuum. The heat flux due to radiation can also occur between  two points of a single body. See \citet{Arpaci1966_01,Reddy1993_01,Nikishkov2010_01,Lienhard2013_01} and \citet{Nithiarasu2016_01} for more discussion.
\subsubsection{Heat transfer}
Heat transfer between a surface and its environment $\bqd_{\text{h}}$ can be written as \citep{Lienhard2013_01}
\eqb{lll}
\ds \bqd_{\text{h}}\cdot\bnu \dis \ds -h\left(T-T_{\text{ref}}\right)~,
\eqe
where $h$ is the heat transfer coefficient and $T_\text{ref}$ is the environment temperature.
\subsection{Examples of the Helmholtz free energy}
The Helmholtz free energy for 3D continua can be written as
\eqb{lll}
\ds \psid \is \ds \psid_{\text{e}}(\bCd_{\text{e}},T)+\psid_{\text{T}}(T)~.
\eqe
For example, $\psid_{\text{e}}$ can be written as \citep{Wriggers2006_01}
\eqb{lll}
\ds \psid_{\text{e}} \is \ds \frac{1}{\rhoId}\left[\frac{\mud}{2}\,\left(\tr\left(\bCd_{\text{e}}\right)-3\right)-\mud\,\ln \Jed+\frac{\lambdad}{2}\,\left(\ln\Jed\right)^2\right]~,
\eqe
where $\bCd_{\text{e}}$ is the elastic part of the right Cauchy--Green deformation tensor and $\Jed=\detd\,\bCd_{\text{e}}$. $\mud$ and $\lambdad$ are material constants and can be functions of temperature, e.g.~\citep{Boyce1992_01}
\eqb{lll}
\mud \is \ds \mu_{0}\,e^{-c_{2}(T-T_{\text{ref}})}~,
\eqe
and $\psid_{\text{T}}$ can be written as \citep{Holzapfel2000_01,Aldakheel2017_01}
\eqb{lll}
\ds \psid_{\text{T}} \is \ds \frac{c_{1}}{\rhoId}\,\left[(T-T_{0})-T\,\ln\,\left(\frac{T}{T_{0}}\right) \right]~.
\eqe
Here, $T_{\text{ref}}$ is the reference temperature and $\mu_{0}$ is the reference shear modulus, i.e.~the shear modulus at $T_{\text{ref}}$, and $c_{1}$ and $c_{2}$ are material constants.
The second Piola-Kirchhoff and entropy can be computed using Eqs.~\eqref{e:3D_nonpolar_S} and \eqref{e:3D_nonpolar_s}, respectively.
The Helmholtz free energy for 2D continua can be written as
\eqb{lll}
\ds \psiz \is \ds \psiez(\bCez,\bkappa,T)+\psiTz(T)~,
\eqe
with the elastic part\footnote{\citet{Zimmermann2018_01} use the alternative form $c_3\,(\buab-\bIuab)(A^{\alpha\gamma}\,b_{\gamma\delta}\,A^{\delta\beta}-b^{\alpha\beta}_{0})$ for the bending part of the Helmholtz free energy~.}
\eqb{lll}
\ds \psiez \is \ds \frac{1}{\rhoIz}\left[\frac{K}{4}\left(\Jez^2-1-2\ln \Jez\right)+\frac{\mu}{2}\left(\frac{\tr(\bCez)}{\Jez}-2\right)\right]+c_3\,\bkappa:\bkappa~,
\eqe
and thermal part \citep{Aldakheel2017_01}
\eqb{lll}
\ds \psiTz \is \ds \frac{c_{1}}{\rhoIz}\,\left[(T-T_{0})-T\,\ln\,\left(\frac{T}{T_{0}}\right) \right]~,
\eqe
where $\mu$ and $K$ can be a function of the temperature, and $c_{1}$ and $c_{3}$ are material constants. $\bsig^{\sharp\lhd}$, $\mu^{\sharp\lhd}$ and $s$ can be computed from Eqs.~\eqref{e:surf_sig_back}, \eqref{e:surf_mu_back} and \eqref{e:surf_s}, see \ref{s:surf_objects_linearization}.
\section{Conclusion}\label{s:conclusion}
A thermomechanical formulation for polar continua under large deformations is derived. It is based on the multiplicative decomposition of the deformation gradient. The proposed formulation for three dimensional polar continua is simplified to three dimensional non-polar continua and Kirchhoff-Love shells. The shell formulation is developed in a tensorial form that is suitable for numerical implementation in both curvilinear and Cartesian coordinates. All formulations consider anisotropic constitutive laws. The weak forms are derived for three dimensional, non-polar continua and Kirchhoff-Love shells.\\
Several aspects of the linearization and discretization of the weak forms are open issues. Special attention is needed for the temperature discretization through the thickness since it drastically changes the condition number of the problem \citep{Surana1992_01}. The numerical results and condition number of the different splitting operators should be considered and compared with a monolithic approach. Analytical benchmark solutions are needed for a validation of the numerical implementation.

\section*{Acknowledgement}{Financial support from the German Research Foundation (DFG) through grant GSC 111 is
gratefully acknowledged. The authors would also like to thank Mr.~Amin Rahmati for helpful discussions.}

\appendix

\section{Strain measures}\label{s:strain_measure}
The additive decomposition of finite strain is only possible for a logarithmic strain. Otherwise the multiplicative decomposition of the deformation gradient is needed for finite strains.
The generalized Lagrangian and Eulerian strain can be written as \citep{Seth1961_01}
\eqb{lll}
\bEd^{(n)} \is \ds \frac{1}{n}\left[\left(\bFd^{\text{T}}\,\bFd\right)^{\frac{n}{2}}-\bolds{1}\right]=\frac{1}{n}\left[\bUd^{n}-\bolds{1}\right]~
\eqe
and
\eqb{lll}
\bed^{(n)} \is \ds \frac{1}{n}\left[\left(\bFd\,\bFd^{\text{T}}\right)^{\frac{n}{2}}-\bolds{1}\right]=\frac{1}{n}\left[\bVd^{n}-\bolds{1}\right]~.
\eqe
The Green-Lagrange $\bE^{(2)}$ and Almansi strain $\be^{(-2)}$ are
\eqb{lll}
\bEd^{(2)} \is \ds \frac{1}{2}\left[\left(\bFd^{\text{T}}\,\bFd\right)-\bolds{1}\right]=\frac{1}{2}\left[\bUd^{2}-\bolds{1}\right]~
\eqe
and
\eqb{lll}
\bed^{(-2)} \is \ds \frac{1}{2}\left[\bolds{1}-\left(\bFd\,\bFd^{\text{T}}\right)^{-1}\right]=\frac{1}{2}\left[\bolds{1}-\bVd^{-2}\right]~.
\eqe
The Lagrangian and Eulerian Hencky strain can be obtained for $n=0$ as
\eqb{lll}
\bEd^{(0)} \is \ds \lim_{n\to 0} \frac{1}{n}\left[\exp{(n\,\ln \bUd)}-\bolds{1}\right]=\ln(\bUd)~
\eqe
and
\eqb{lll}
\bed^{(0)} \is \ds \lim_{n\to 0} \frac{1}{n}\left[\exp{(n\,\ln \bVd)}-\bolds{1}\right]=\ln(\bVd)~
\eqe
where L'H{\^o}pital's rule is used. For two consecutive deformation gradients of $\bFd_{1}$ and $\bFd_{2}$ such that $\bFd=\bFd_{2}\,\bFd_{1}$, $\bEd^{(n)}=\bEd^{(n)}_{1}+\bEd^{(n)}_{2}$ is true only if the logarithmic strain $n=0$ is used. Here, $\bE^{(n)}_{i}$ is defined as
\eqb{lll}
\bEd^{(n)}_{i} \dis \ds \frac{1}{n} \left[\left(\bFd^{\text{T}}_{i}\,\bFd_{i}\right)^{\frac{n}{2}}-\bolds{1}\right]~.
\eqe
It means that the logarithmic strain can be additively decomposed for finite deformations \citep{Bruhns2015_01}. $\bE^{(n)}$, for $n\neq0$, can be written as
\eqb{lll}
\ds \bEd^{(n)} \is \ds \left(\bFd^{\text{T}}_{1}\right)^{\frac{n}{2}}\,\bEd^{(n)}_{2}\,\left(\bFd_{1}\right)^{\frac{n}{2}}+\bEd^{(n)}_{1}~.
\eqe
A similar relation can be obtained if $\bFd_{1}$ and $\bFd_{2}$ are assumed to be the thermal and elastic deformation gradients. The Green-Lagrange strain can be decomposed additively if $\bFd_{1}\simeq\bolds{1}$. The logarithmic strain facilitates the development of material models, but the numerical implementation is complicated \citep{Ghaffari2018_03}. A similar discussion can be provided for the decomposition of surface strain measures. This multiplicative formulation can be generalized to thermoplasticity as \citep{Meissonnier2001_01}
\eqb{lll}
\bFd \is \bFd_{\text{e}}\,\bFd_{\text{p}}\,\bFd_{\text{T}}~,
\eqe
where $\bFd_{\text{p}}$ is the plastic deformation gradient\footnote{The alternative decomposition $\bFd=\bFd_{\text{e}}\,\bFd_{\text{T}}\,\bFd_{\text{p}}$ can also be used \citep{McHugh1993_01,Srikanth1999_01}. See \citet{Wang2017_01} for including the phase transformation deformation gradient $\bFd_{\text{pt}}$ as $\bFd=\bFd_{\text{e}}\,\bFd_{\text{pt}}\,\bFd_{\text{T}}$.}.

\section{Linearization of surface objects}\label{s:surf_objects_linearization}
Here, the tensorial derivatives of surface objects are obtained, so the componentwise derivatives of \citet{Sauer2017_01} is avoided. Hence, the linearization can be used in a curvilinear or Cartesian system. The derivative of $\Jz$ and $\bCz^{-1}$ with respect to $\bCz$ can be written as
\eqb{lll}
\ds \pa{\Jz}{\bCz} \is \ds \pa{\sqrt{\det{\bCz}}}{\bCz} =\frac{\Jz}{2}\bCz^{-\text{T}}=\frac{\Jz}{2}\bCz^{-1}~,
\eqe
and
\eqb{lll}
\ds \frac{\partial\bCz^{-1}}{\oplus\partial\bCz}\is \ds -\frac{1}{2}\left(\bCz^{-1}\otimes\bCz^{-1}+\bCz^{-1}\boxtimes\bCz^{-1}\right)~.
\eqe
To obtain the derivative of $H$, it is written as
\eqb{lll}
\ds H \is \ds \frac{1}{2}\tr(\bb)=\frac{1}{2}\bb:\bi=\frac{1}{2}\buab\,\vaa\otimes\vab:\agd\,\vaug\otimes\vaud =\frac{1}{2}\bb^{\flat\lhd}:\bCz^{-1}~,
\eqe
with
\eqb{lll}
\ds \bb^{\flat\lhd} \is \buab\,\vAa\otimes\vAb~.
\eqe
The derivatives of $H$ with respect to $\bCz$ and $\bb^{\flat\lhd}$ are
\eqb{lll}
\ds\pa{H}{\bCz} \is \ds \frac{1}{2} \bb^{\flat\lhd}\bullet\circ\left[-\frac{1}{2}\left(\bCz^{-1}\otimes\bCz^{-1}+\bCz^{-1}\boxtimes\bCz^{-1}\right)\right] =-\frac{1}{2}\bb^{\sharp\lhd}~,
\eqe
\eqb{lll}
\ds \pa{H}{\bb^{\flat\lhd}} \is \ds \frac{1}{2}\bCz^{-1}~,
\eqe
with
\eqb{lll}
\bb^{\sharp\lhd} \is \bab\,\vAua\otimes\vAub=\bCz^{-1}\,\bb^{\flat\lhd}\,\bCz^{-1}~.
\eqe
Similarly, $\kappa$ can be written as
\eqb{lll}
\kappa \is \det(\bb)=\det\big(\buab\,\vaa\otimes\vab\big)=\det\left(\buab\,\bFz^{-\text{T}}\vAa\otimes\vAb\,\bFz^{-1}\right)= \det\left(\bFz^{-\text{T}}\,\bb^{\flat\lhd}\,\bFz^{-1}\right)\\[4mm]
\is \det\big(\bb^{\flat\lhd}\big)\det\big(\bCz^{-1}\big)~,
\eqe
and its derivatives with respect to $\bCz$ and $\bb^{\flat\lhd}$ are
\eqb{lll}
\ds\pa{\kappa}{\bCz} =\kappa_{|\bCz} \is \ds -\kappa\,\bCz^{-1}~,
\eqe
\eqb{lll}
\ds\pa{\kappa}{\bb^{\flat\lhd}}= \kappa_{|\bb^{\flat\lhd}} \is \ds \kappa\,\bb^{\flat\lhd-1}~.
\eqe
Finally, the derivatives of $\bb^{\sharp\lhd}$ with respect to $\bCz$ and $\bb^{\flat\lhd}$ are
\eqb{lll}
\ds  \frac{\partial\bb^{\sharp\lhd}}{\oplus\partial\bCz}\is \ds -\frac{1}{2}\left(\bCz^{-1}\otimes\bb^{\sharp\lhd}+\bCz^{-1}\boxtimes\bb^{\sharp\lhd} +\bb^{\sharp\lhd}\otimes\bCz^{-1}+\bb^{\sharp\lhd}\boxtimes\bCz^{-1}\right)~,
\eqe
\eqb{lll}
\ds \frac{\partial\bb^{\sharp\lhd}}{\oplus\partial\bb^{\flat\lhd}}\is \ds \frac{1}{2}\left(\bCz^{-1}\otimes\bCz^{-1}+\bCz^{-1}\boxtimes\bCz^{-1}\right)~.
\eqe
\section{Time derivatives of surfaces objects}
The material time derivatives of the tangent vectors can be written as (see \citet[p.~657-681]{Bower2009_01}, \citet{Sahu2017_01} or \citet{Wriggers2006_01})
\eqb{lll}
\ds \frac{\text{D} \ba_{\alpha}}{\text{D} t} \is \ds \pa{ \bv}{\xi^{\alpha}}=\pa{v^{i}}{\xi^{\alpha}}+v^{i}\pa{\bg_{i}}{\xi^{\alpha}} =\pa{v^{i}}{\xi^{\alpha}}\,\bg_{i}+v^{i}\Gamma^{k}_{i\alpha}\bg_{k}\\[4mm]
\is w_{\alpha}^{~\beta}\,\vaub+w_{\alpha}\bn\\[4mm]
\is w_{\alpha\beta}\,\vab+w_{\alpha}\bn~,
\eqe
with
\eqb{lll}
w_{\alpha}^{~\beta} \is v^{\beta}_{;\alpha}-v\,b^{\beta}_{\alpha}~,
\eqe
\eqb{lll}
w_{\alpha} \is v^{\lambda}\,b_{\lambda\alpha}+v_{,\alpha}~,
\label{e:surf_w_a}
\eqe
\eqb{lll}
w_{\alpha\beta} \is w_{\alpha}^{~\mu}\,a_{\mu\beta}~.
\label{e:surf_w_ab}
\eqe
Taking the time derivative of $\vaua\cdot\bn=0$ and $\bn\cdot\bn$ results in
\eqb{lll}
\dot{\bn}\cdot\vaua=-\dot{\ba}_{\alpha}\cdot\bn~,
\eqe
\eqb{lll}
\dot{\bn}\cdot\bn=0~,
\eqe
where $\dot{\bn}\cdot\bn=0$ indicates that $\dot{\bn}$ has no out of plane components, so
\eqb{lll}
\dot{\bn}=-\ba^{\alpha}\otimes\bn\cdot\dot{\ba}_{\alpha}=-w_{\alpha}\,\vaa=-w^{\alpha}\,\vaua~,
\eqe
\eqb{lll}
w^{\alpha} \is w_{\beta}\,a^{\beta\alpha} ~.
\eqe
The velocity gradient for the mid-surface $\blIz$ can be written as
\eqb{lll}
\ds \blIz = \ds \dot{\bFIz}\,{\bFIz}^{-1} \is \ds \dot{\ba}_{\alpha}\otimes\ba^{\alpha}+\frac{\dot{\lambda}_{3}}{\lambda_{3}}\,\bn\otimes\bn+\,\dot{\bn}\otimes\bn\\[3mm]
\is \ds w_{\alpha\beta}\,\vab\otimes\vaa+w_{\alpha}\bn\otimes\vaa+\,\dot{\bn}\otimes\bn +\frac{\dot{\lambda}_{3}}{\lambda_3}\,\bn\otimes\bn\\[3mm]
\is \ds w_{\alpha\beta}\,\vab\otimes\vaa-\bn\otimes\dot{\bn}+\,\dot{\bn}\otimes\bn +\frac{\dot{\lambda}_{3}}{\lambda_3}\,\bn\otimes\bn~.
\label{e:surf_velocity_grad}
\eqe
Furthermore, $\dot{a}_{\alpha\beta}$ and $\dot{b}_{\alpha\beta}$ can be written as
\eqb{lll}
\dot{a}_{\alpha\beta} \is  \dot{\ba}_{\alpha}\cdot\vaub+\vaua\cdot\dot{\ba}_{\beta}=w_{\alpha\beta}+w_{\beta\alpha}~,
\label{e:auab_rate}
\eqe
\eqb{lll}
\dot{b}_{\alpha\beta} \is \ds \dot{\ba}_{\alpha;\beta}\cdot\bn+\ba_{\alpha;\beta}\cdot\dot{\bn}\\[3mm] \is w_{\alpha}^{~\gamma}\,b_{\gamma\beta}+w_{\alpha;\beta}=w_{\alpha\gamma}\,b^{\gamma}_{\beta}+w_{\alpha;\beta}~.
\label{e:buab_rate}
\eqe

\bibliographystyle{model1-num-names}
\bibliography{bibliography}

\end{document}

%% file: neco_updated.tex
\newcommand {\detd}{{\det}}
\newcommand {\curld}{{\text{curl}}}
\newcommand {\divd}{{\text{div}}}
\newcommand {\Divd}{{\text{Div}}}
\newcommand {\gradd}{{\text{grad}}}
\newcommand {\Gradd}{{\text{Grad}}}

\newcommand {\detz}{\text{det}_{\text{s}}}
\newcommand {\divz}{\text{div}_{\text{s}}}
\newcommand {\Divz}{\text{Div}_{\text{s}}}
\newcommand {\gradz}{\text{grad}_{\text{s}}}
\newcommand {\Gradz}{\text{Grad}_{\text{s}}}

\newcommand {\bad}{{\ba}}
\newcommand {\bAd}{{\bA}}
\newcommand {\balphad}{{\bolds{\alpha}}}
\newcommand {\bbd}{{\bb}}
\newcommand {\bBd}{{\bB}}
\newcommand {\bcd}{{\bc}}
\newcommand {\bCd}{{\bC}}
\newcommand {\bCTd}{{\bC}_{\text{T}}}
\newcommand {\bCed}{{\bC}_{\text{e}}}
\newcommand {\GammaId}[3]{{\Gamma}_{0\,#1#2}^{#3}}
\newcommand {\Gammad}[3]{{\Gamma}_{#1#2}^{#3}}
\newcommand {\bdd}{{\bd}}
\newcommand {\bed}{{\be}}
\newcommand {\bEd}{{\bE}}
\newcommand {\bsEd}{{\bolds{\sE}}}
\newcommand {\bffd}{{\bff}}
\newcommand {\bFd}{{\bF}}
\newcommand {\bFTd}{{\bF}_{\!\text{T}}}
\newcommand {\bFed}{{\bF}_{\!\text{e}}}
\newcommand {\bFTTd}{{\bF}_{\!\text{T},T}}
\newcommand {\bgTdui}[1]{\gbg_{#1}}
\newcommand {\bgTdi}[1]{\gbg^{#1}}
\newcommand {\gTdui}[1]{\ggn_{#1}}
\newcommand {\bHd}{{\bH}}
\newcommand {\Jd}{{J}}
\newcommand {\JTd}{{J}_{\text{T}}}
\newcommand {\Jed}{{J}_{\text{e}}}
\newcommand {\bkd}{{\bk}}
\newcommand {\lambdad}{{\lambda}}
\newcommand {\bld}{{\bl}}
\newcommand {\blTd}{{\bl}_{\text{T}}}
\newcommand {\bled}{{\bl}_{\text{e}}}
\newcommand {\bLd}{{\bL}}
\newcommand {\bmd}{{\bm}}
\newcommand {\mud}{{\mu}}
\newcommand {\bmurd}{{\bar{\bmu}}}
\newcommand {\bmurId}{{\bar{\bmu}}_{0}}
\newcommand {\bPd}{{\bP}}
\newcommand {\rhod}{{\rho}}
\newcommand {\rhoId}{{\rho}_{0}}
\newcommand {\bsigd}{{\bsig}}
\newcommand {\bqd}{{\bq}}
\newcommand {\bQd}{{\bQ}}
\newcommand {\rd}{{r}}
\newcommand {\sd}{{s}}
\newcommand {\bSd}{{\bS}}
\newcommand {\bSTd}{{\bS}_{\text{T}}}
\newcommand {\btd}{{\bt}}
\newcommand {\btId}{{\bt}_{0}}
\newcommand {\ud}{{u}}
\newcommand {\bud}{{\bu}}
\newcommand {\bUd}{{\bU}}
\newcommand {\bvd}{{\bv}}
\newcommand {\bVd}{{\bV}}

\newcommand {\bvarphid}{{\bolds{\varphi}}}
\newcommand {\bwd}{{\bw}}
\newcommand {\bwAz}{\bw^{\bowtie}}
\newcommand {\bwAd}{{\bw}^{\bowtie}}
\newcommand {\bxd}{{\bx}}
\newcommand {\byd}{{\by}}
\newcommand {\bXd}{{\bX}}
\newcommand {\bXTd}{\gbX}
\newcommand {\psid}{{\psi}}
\newcommand {\Psid}{{\Psi}}

\newcommand{\bbCd}{{\mathbb{C}}}

\newcommand {\aTz}{\mathfrak{a}}
\newcommand {\bTzui}[1]{\mathfrak{b}_{#1}}
\newcommand {\baTz}{\boldsymbol{\mathfrak{a}}}
\newcommand {\bbTz}{\boldsymbol{\mathfrak{b}}}
\newcommand {\baTzui}[1]{\gba_{#1}}
\newcommand {\baTzi}[1]{\gba^{#1}}
\newcommand {\bcz}{\bc_{\text{s}}}
\newcommand {\bCz}{\bC_{\text{s}}}
\newcommand {\bCzt}{\dot{\bC}_{\text{s}}}
\newcommand {\bCTz}{\bC_{\text{sT}}}
\newcommand {\bCez}{\bC_{\text{se}}}
\newcommand {\bCezt}{\dot{\bC}_{\text{se}}}
\newcommand {\bCaz}{\accentset{\ast}{\bC}}
\newcommand {\bCaez}{\accentset{\ast}{\bC}_{\!\text{e}}}
\newcommand {\bCaTz}{\accentset{\ast}{\bC}_{\!\text{T}}}
\newcommand {\bCcaz}{\accentset{\circledast}{\bC}}
\newcommand {\bCcaez}{\accentset{\circledast}{\bC}_{\!\text{e}}}
\newcommand {\bCcaTz}{\accentset{\circledast}{\bC}_{\!\text{T}}}
\newcommand {\bffz}{\bff_{\!\text{s}}}
\newcommand {\bFz}{\bF_{\!\text{s}}}
\newcommand {\bFTz}{\bF_{\!\text{sT}}}
\newcommand {\bFez}{\bF_{\!\text{se}}}
\newcommand {\bFaz}{\accentset{\ast}{\bF}}
\newcommand {\bFaTz}{\accentset{\ast}{\bF}_{\!\text{T}}}
\newcommand {\bFaez}{\accentset{\ast}{\bF}_{\!\text{e}}}
\newcommand {\bFcaz}{\accentset{\circledast}{\bF}}
\newcommand {\bFcaTz}{\accentset{\circledast}{\bF}_{\!\text{T}}}
\newcommand {\bFcaez}{\accentset{\circledast}{\bF}_{\!\text{e}}}
\newcommand {\bFIz}{\accentset{0}{\bF}}
\newcommand {\bFTzi}[1]{\bF_{\!\text{sT}#1}}
\newcommand {\bFTTz}{\bF_{\!\text{sT},T}}

\newcommand {\gcaz}{\accentset{\circledast}{g}}
\newcommand {\bGcaz}{\accentset{\circledast}{\bG}}
\newcommand {\bgcaz}{\accentset{\circledast}{\bg}}
\newcommand {\bgcaTz}{\accentset{\circledast}{\gbg}}
\newcommand {\bGcazF}{\overset{\circledast}{\bG}}
\newcommand {\bgcazF}{\overset{\circledast}{\bg}}
\newcommand {\bgcaTzF}{\overset{\circledast}{\gbg}}
\newcommand {\bgcaTzui}[1]{\accentset{\circledast}{\gbg}_{#1}}

\newcommand {\gcaTzi}[1]{\accentset{\circledast}{\ggn}_{#1}}
\newcommand {\Gcaz}{\accentset{\circledast}{G}}
\newcommand {\GammaIz}[3]{\Gamma_{\text{s}0\,#1#2}^{#3}}
\newcommand {\Gammaz}[3]{\Gamma_{\text{s}\,#1#2}^{#3}}
\newcommand {\bHz}{\bH_{\text{s}}}
\newcommand {\Jz}{J_{\text{s}}}
\newcommand {\JTz}{J_{\text{sT}}}
\newcommand {\Jez}{J_{\text{se}}}
\newcommand {\bkz}{\bk_{\text{s}}}
\newcommand {\blz}{\bl_{\text{s}}}
\newcommand {\blTz}{\bl_{\text{sT}}}
\newcommand {\blez}{\bl_{\text{se}}}
\newcommand {\blIz}{\accentset{0}{\bl}}
\newcommand {\blcaz}{\accentset{\circledast}{\bl}}

\newcommand {\murz}{\bar{\mu}}
\newcommand {\bmz}{\bm_{\text{s}}}
\newcommand {\bmuz}{\bmu_{\text{s}}}
\newcommand {\bmuIz}{\bmu_{\text{0s}}}
\newcommand {\bmurz}{\bar{\bmu}_{\text{s}}}
\newcommand {\rhoz}{\rho_{\text{s}}}
\newcommand {\rhoIz}{\rho_{\text{s0}}}
\newcommand {\rz}{r_{\text{s}}}
\newcommand {\bTz}{\bT_{\text{s}}}
\newcommand {\bxz}{\bx_{\text{s}}}
\newcommand {\bxzi}[1]{\bx_{\text{s}#1}}
\newcommand {\bxcaz}{\accentset{\circledast}{\bx}}
\newcommand {\bXz}{\bX_{\!\text{s}}}
\newcommand {\bXTz}{\gbX_{\text{s}}}
\newcommand {\bXcaz}{\accentset{\circledast}{\bX}}
\newcommand {\bXcaTzi}[1]{\accentset{\circledast}{\gbX}_{#1}}
\newcommand {\bXcaTz}{\accentset{\circledast}{\gbX}}
\newcommand {\Psiz}{\Psi_{\text{s}}}
\newcommand {\psiz}{\psi_{\text{s}}}
\newcommand {\psizt}{\dot{\psi}_{\text{s}}}
\newcommand {\psiTz}{\psi_{\text{sT}}}
\newcommand {\psiez}{\psi_{\text{se}}}
\newcommand {\bqz}{\bq_{\text{s}}}
\newcommand {\bqcaz}{\accentset{\circledast}{\bq}}
\newcommand {\sz}{s_{\text{s}}}
\newcommand {\szt}{\dot{s}_{\text{s}}}
\newcommand {\bSz}{\bS_{\text{s}}}
\newcommand {\bsigz}{\bsig_{\text{s}}}
\newcommand {\btz}{\bt_{\text{s}}}
\newcommand {\btIz}{\bt_{\text{s}0}}
\newcommand {\uz}{u_{\text{s}}}
\newcommand {\uzt}{\dot{u}_{\text{s}}}
\newcommand {\buz}{\bu_{\text{s}}}
\newcommand {\buzi}[1]{\bu_{\text{s}#1}}
\newcommand {\bUz}{\bU_{\!\text{s}}}
\newcommand {\bUzi}[1]{\bU_{\!\text{s}#1}}
\newcommand {\bvz}{\bv_{\text{s}}}
\newcommand {\bvzt}{\dot{\bv}_{\text{s}}}
\newcommand {\bvzi}[1]{\bv_{\!\text{s}#1}}
\newcommand {\bVz}{\bV_{\!\text{s}}}
\newcommand {\bVzi}[1]{\bV_{\!\text{s}#1}}

\newcommand{\bbCz}{\mathbb{C}_{\text{s}}}
\newcommand{\mcalScaz}{\accentset{\circledast}{\mcalS}}

\newcommand{\mcalScazF}{\overset{\circledast}{\mcalS}}
\newcommand{\mcalScaTzF}{\overset{\circledast}{\mcalS}_{\text{T}}}
\newcommand{\mcalSTz}{\mcalS_\text{T}}
\newcommand{\mcalBTd}{\mcalB_\text{T}}

\newcommand {\ggn}{\mathfrak{g}}

\newcommand {\gba}{\boldsymbol{\mathfrak{a}}}

\newcommand {\gbg}{\boldsymbol{\mathfrak{g}}}

\newcommand {\gbX}{\boldsymbol{\mathfrak{X}}}

\newcommand{\bitm}{\begin{itemize}}
\newcommand{\eitm}{\end{itemize}}
\newcommand{\bnumr}{\begin{enumerate}}
\newcommand{\enumr}{\end{enumerate}}
\newcommand{\bolds}[1]{\boldsymbol{#1}}

\newcommand {\sigab}{\sigma^{\alpha\beta}}

\newcommand {\aab}{a^{\alpha\beta}}
\newcommand {\auab}{a_{\alpha\beta}}
\newcommand {\agd}{a^{\gamma\delta}}

\newcommand {\abg}{a^{\beta\gamma}}

\newcommand {\aag}{a^{\alpha\gamma}}

\newcommand {\augd}{a_{\gamma\delta}}

\newcommand {\adb}{a^{\delta\beta}}

\newcommand {\Aab}{A^{\alpha\beta}}
\newcommand {\Auab}{A_{\alpha\beta}}
\newcommand {\Aag}{A^{\alpha\gamma}}

\newcommand {\Abg}{A^{\beta\gamma}}

\newcommand {\Augd}{A_{\gamma\delta}}

\newcommand {\Adb}{A^{\delta\beta}}

\newcommand {\Guij}{G_{ij}}
\newcommand {\Gij}{G^{ij}}

\newcommand {\guij}{g_{ij}}
\newcommand {\gij}{g^{ij}}

\newcommand {\Mab}{M^{\alpha\beta}}
\newcommand {\Nab}{N^{\alpha\beta}}

\newcommand {\bab}{b^{\alpha\beta}}

\newcommand {\buab}{b_{\alpha\beta}}

\newcommand {\bIuab}{b_{0\,\alpha\beta}}

\newcommand {\tauab}{\tau^{\alpha\beta}}

\newcommand{\mcalB}{\mathcal{B}}

\newcommand{\mcalS}{\mathcal{S}}

\newcommand {\eqb}[1]{\begin{equation}\begin{array}{#1}}
\newcommand {\eqe}{\end{array}\end{equation}}

\newcommand {\esb}[1]{\begin{equation*}\begin{array}{#1}}
\newcommand {\ese}{\end{array}\end{equation*}}
\newcommand {\ds}{\displaystyle}

\newcommand {\pa}[2]{\frac{\partial{#1}}{\partial{#2}}}

\newcommand {\paqq}[3]{\frac{\partial^2{#1}}{\partial{#2}\,\partial{#3}}}
\newcommand {\back}{\! \! \!}
\newcommand {\is}{\back &=& \back}
\newcommand {\dis}{\back &:=& \back}

\newcommand {\plus}{\back &+& \back}

\newcommand {\tr}{\mathrm{tr}\,}

\newcommand {\dif}{\mathrm{d}}

\newcommand {\II}{{I\kern-.3em I}}
\newcommand {\III}{{I\kern-.3em I\kern-.3em I}}

\newcommand {\ba}{\boldsymbol{a}}
\newcommand {\bb}{\boldsymbol{b}}
\newcommand {\bc}{\boldsymbol{c}}
\newcommand {\bd}{\boldsymbol{d}}
\newcommand {\be}{\boldsymbol{e}}
\newcommand {\bff}{\boldsymbol{f}}
\newcommand {\bg}{\boldsymbol{g}}

\newcommand {\bi}{\boldsymbol{i}}

\newcommand {\bk}{\boldsymbol{k}}
\newcommand {\bl}{\boldsymbol{l}}

\newcommand {\bm}{\boldsymbol{m}}
\newcommand {\bn}{\boldsymbol{n}}

\newcommand {\bq}{\boldsymbol{q}}

\newcommand {\bt}{\boldsymbol{t}}
\newcommand {\bu}{\boldsymbol{u}}
\newcommand {\bv}{\boldsymbol{v}}
\newcommand {\bw}{\boldsymbol{w}}
\newcommand {\bx}{\boldsymbol{x}}
\newcommand {\by}{\boldsymbol{y}}

\newcommand {\bmu}{\boldsymbol{\mu}}

\newcommand {\bnu}{\mbox{\boldmath$\nu$}}

\newcommand {\bA}{\boldsymbol{A}}
\newcommand {\bB}{\boldsymbol{B}}
\newcommand {\bC}{\boldsymbol{C}}

\newcommand {\bE}{\boldsymbol{E}}
\newcommand {\bF}{\boldsymbol{F}}
\newcommand {\bG}{\boldsymbol{G}}
\newcommand {\bH}{\boldsymbol{H}}
\newcommand {\bI}{\boldsymbol{I}}

\newcommand {\bL}{\boldsymbol{L}}

\newcommand {\bN}{\boldsymbol{N}}

\newcommand {\bP}{\boldsymbol{P}}
\newcommand {\bQ}{\boldsymbol{Q}}

\newcommand {\bS}{\boldsymbol{S}}
\newcommand {\bT}{\boldsymbol{T}}
\newcommand {\bU}{\boldsymbol{U}}
\newcommand {\bV}{\boldsymbol{V}}

\newcommand {\bX}{\boldsymbol{X}}

\newcommand {\bsig}{\mbox{\boldmath$\sigma$}}

\newcommand {\btau}{\mbox{\boldmath$\tau$}}

\newcommand {\bone}{\mathbf{1}}

\newcommand {\bzero}{\mathbf{0}}

\newcommand {\bbR}{\mathbb{R}}

\newcommand {\IR}{{\rm\kern.24em
   \vrule width.02em height1.53ex depth-.05ex
   \kern-.3em R}}
\newcommand {\ic}{{\rm\kern.20em
   \vrule width.02em height1.0ex depth-.05ex
   \kern-.22em c}}
\newcommand {\ia}{{\rm\kern.20em
   \vrule width.02em height1.05ex depth-.0ex
   \kern-.25em a}}
\newcommand {\IC}{{\rm\kern.24em
   \vrule width.02em height1.4ex depth-.05ex
   \kern-.26em C}}
\newcommand {\ID}{{\rm\kern.34em
   \vrule width.02em height1.5ex depth-.05ex
   \kern-.36em D}}
\newcommand {\IS}{{\rm\kern.24em
   \vrule width.02em height1.6ex depth.05ex
   \kern-.26em S}}
\newcommand {\IT}{{\rm\kern.50em
   \vrule width.02em height1.55ex depth-.05ex
   \kern-.52em T}}

\newcommand {\IE}{{\rm\kern.24em
   \vrule width.02em height1.55ex depth-.05ex
   \kern-.33em E}}
\newcommand {\IEa}{{\rm\kern.24em
   \vrule width.02em height1.55ex depth-.05ex
   \kern-.33em E}^{1}_{ijkl}}
\newcommand {\IEb}{{\rm\kern.24em
   \vrule width.02em height1.55ex depth-.05ex
   \kern-.33em E}^{2}_{ijkl}}

\newcommand {\sB}{\mathcal{B}}

\newcommand {\sE}{\mathcal{E}}

\newcommand {\sK}{\mathcal{K}}

\newcommand {\sP}{\mathcal{P}}
\newcommand {\sQ}{\mathcal{Q}}

\newcommand {\sS}{\mathcal{S}}

\newcommand {\sU}{\mathcal{U}}
\newcommand {\sV}{\mathcal{V}}

\newcommand {\vaua}{\ba_\alpha}
\newcommand {\vaa}{\ba^\alpha}

\newcommand {\vaub}{\ba_\beta}
\newcommand {\vab}{\ba^\beta}

\newcommand {\vaug}{\ba_\gamma}
\newcommand {\vag}{\ba^\gamma}

\newcommand {\vaud}{\ba_\delta}

\newcommand {\vAua}{\bA_\alpha}
\newcommand {\vAa}{\bA^\alpha}

\newcommand {\vAub}{\bA_\beta}
\newcommand {\vAb}{\bA^\beta}

\newcommand {\vAug}{\bA_\gamma}
\newcommand {\vAg}{\bA^\gamma}

\newcommand {\Ass}[2]{\kern 0.9ex \vrule width0.45em height0.2ex depth0ex \kern -2.1ex \bigwedge_{#1}^{#2}}
\newcommand {\ASS}[2]{\kern 1.45ex \vrule width0.5em height0.2ex depth0ex \kern -2.65ex \bigwedge_{#1}^{#2}}

\newcommand {\bkappa}{\boldsymbol{\kappa}}

